\documentclass [12pt]{article}
\usepackage{latexsym,amssymb,amsmath,mathrsfs}
\usepackage{pst-all}
\usepackage{graphicx}

\newtheorem{thm}{Theorem}[section]
\newtheorem{prop}[thm]{Proposition}
\newtheorem{cor}[thm]{Corollary}
\newtheorem{lem}[thm]{Lemma}
\newtheorem{defn}[thm]{Definition}
\newtheorem{remark}[thm]{Remark}
\newtheorem{example}[thm]{Example}

\makeatletter \@addtoreset{equation}{section} \makeatother

\newcommand{\nn}{\nonumber}
\newcommand{\ds}{\displaystyle}

\newcommand{\bthm}[1]{\begin{thm}\label{th:#1}}
\newcommand{\bprop}[1]{\begin{prop}\label{prop:#1}}
\newcommand{\blem}[1]{\begin{lem}\label{lem:#1}}
\newcommand{\bcor}[1]{\begin{cor}\label{cor:#1}}
\newcommand{\bdefn}[1]{\begin{defn}\label{def:#1}\rm}
\newcommand{\bdefnn}[2]{\begin{defn}{\bf (#2)}\label{def:#1}\rm}
\newcommand{\brem}[1]{\begin{remark}\label{rem:#1}\rm}
\newcommand{\bexam}[1]{\begin{example}\label{ex:#1}\rm}
\newcommand{\ethm}{\end{thm}}
\newcommand{\eprop}{\end{prop}}
\newcommand{\elem}{\end{lem}}
\newcommand{\ecor}{\end{cor}}
\newcommand{\edefn}{\end{defn}}
\newcommand{\erem}{\end{remark}}
\newcommand{\eexam}{\end{example}}

\newenvironment{demo}[1]{\par\begin{trivlist}%
\item[]{\bf #1}\ }{\end{trivlist}\par}
\newcommand{\bpf}{\begin{demo}{\bf Proof.\ }}
\newcommand{\apf}[1]{\begin{demo}{\bf Proof of #1.\ }}
\newcommand{\epf}{\hfill $\square$ \end{demo}}

\newcommand{\Thm}[1]{Theorem~\ref{th:#1}}
\newcommand{\Prop}[1]{Proposition~\ref{prop:#1}}
\newcommand{\Lem}[1]{Lemma~\ref{lem:#1}}
\newcommand{\Cor}[1]{Corollary~\ref{cor:#1}}
\newcommand{\Defn}[1]{Definition~\ref{def:#1}}
\newcommand{\Rem}[1]{Remark~\ref{rem:#1}}
\newcommand{\Exam}[1]{Example~\ref{ex:#1}}
\newcommand{\eq}[1]{\eqref{eq:#1}}

\renewcommand{\P}{\mathbb{P}}
\newcommand{\E}{\mathbb{E}}
\newcommand{\R}{\mathbb{R}}
\newcommand{\N}{\mathbb{N}}
\newcommand{\sB}{\mathscr{B}}
\newcommand{\sC}{\mathscr{C}}
\newcommand{\sE}{\mathscr{E}}
\newcommand{\sF}{\mathscr{F}}
\newcommand{\sO}{\mathscr{O}}

\newcommand{\gm}{\gamma}
\newcommand{\dl}{\delta}
\newcommand{\ep}{\varepsilon}
\newcommand{\te}{\theta}
\newcommand{\h}{\eta}
\newcommand{\lm}{\lambda}
\newcommand{\ro}{\rho}
\newcommand{\sg}{\sigma}
\newcommand{\ph}{\varphi}

\newcommand{\Ph}{\Phi}
\newcommand{\Om}{\Omega}

\newcommand{\abra}[1]{\left( #1 \right)}
\newcommand{\bbra}[1]{\left\{ #1 \right\}}
\newcommand{\cbra}[1]{\left[ #1 \right]}

\newcommand{\wg}{\wedge}

\newcommand{\vnorm}[1]{\left\|#1\right\|_{\mathrm{var}}}
\newcommand{\varm}[3]
{\vnorm{\P_{#1} \circ Z (#3)^{-1} - \P_{#2} \circ Z (#3)^{-1}}}

\newcommand{\pssierp}{
  \multido{\i=0+4}{3}{
    \rput{\i}{\pspolygon(1;9)(2;11)(1;1)}
  }
}
\newcommand{\pssierpp}{
  \multido{\ia=11+4}{3}{
    \rput(1;\ia){
      \scalebox{0.5}{
        \pssierp
      }
    }
  }
}
\newcommand{\pssierppp}{
  \multido{\ia=11+4}{3}{
    \rput(1;\ia){
      \scalebox{0.5}{
        \pssierpp
      }
    }
  }
}
\newcommand{\pssierpppp}{
  \multido{\ia=11+4}{3}{
    \rput(1;\ia){
      \scalebox{0.5}{
        \pssierppp
      }
    }
  }
}
\newcommand{\pssierppppp}{
  \multido{\ia=11+4}{3}{
    \rput(1;\ia){
      \scalebox{0.5}{
        \pssierpppp
      }
    }
  }
}
\newcommand{\pssierpinski}{
  \multido{\ia=11+4}{3}{
    \rput(1;\ia){
      \scalebox{0.5}{
        \pssierppppp
      }
    }
  }
}

\title{On uniqueness of maximal coupling 
for diffusion processes with a reflection
}
\author{Kazumasa Kuwada
\footnote{
Department of Mathematics, 
Faculty of Science, 
Ochanomizu University, 
Tokyo 112-8610, Japan
}
\footnote{
\textit{Tel}:\texttt{+81-3-5978-5300}, 
\textit{e-mail}: \texttt{kkuwada@math.ocha.ac.jp}
}
}
\date{}
\begin{document}
\maketitle
\raggedbottom 
\begin{abstract}
A maximal coupling of two diffusion processes 
makes two diffusion particles meet as early as possible. 
We study the uniqueness of maximal couplings 
under a sort of `reflection structure' which ensures 
the existence of such couplings. 
In this framework, 
the uniqueness in the class of Markovian couplings 
holds for the Brownian motion on a Riemannian manifold 
whereas it fails in more singular cases. 
We also prove that a Kendall-Cranston coupling is maximal 
under the reflection structure. 
\end{abstract}
\begin{list}{\textbf{Key words:}}{\setlength{\labelwidth}{72pt}
\setlength{\leftmargin}{72pt}}
\item
Diffusion process, 
maximal coupling, 
mirror coupling,
Kendall-Cranston coupling
\end{list}

\section{Introduction} \label{sec:intro} 

The concept of coupling is very useful in various problems in probability. 
Given probability measures $\mu_1$ on $(\Om_1 , \sF_1)$ 
and $\mu_2$ on $(\Om_2 , \sF_2)$, 
we say $\mu$ a coupling of $\mu_1$ and $\mu_2$, 
or $\mu \in \sC ( \mu_1 , \mu_2)$, 
when $\mu$ is a probability measure 
on $(\Om_1 ,\sF_1)  \times ( \Om_2 , \sF_2)$ 
so that its marginal distributions coincide 
with $\mu_1$ and $\mu_2$ respectively. 
That is, $\mu(A_1 \times \Om_2)=\mu_1(A_1)$ for $A_1 \in \sF_1$ and 
$\mu(\Om_1 \times A_2) = \mu_2(A_2)$ for $A_2 \in \sF_2$. 

We consider couplings of a diffusion process 
$( \{ Z(t)\}_{t \ge 0} , \{ \P_x \}_{x \in M} )$ 
on a topological space $X$. 
A coupling $\P \in \sC (\P_x , \P_y )$ determines a stochastic process 
$(Z_1, Z_2)$ on $X \times X$ so that each individual component moves 
as the diffusion process starting at $x$ and $y$ respectively. 
A characteristic of couplings on which we concentrate our attention is 
the coupling time $T(Z_1, Z_2)$, the time when $Z_1$ and $Z_2$ coalesce
(defined in \eq{ctime}). 
In many applications, 
we would like to make the coupling probability $\P[ T > t ]$ small 
by constructing a suitable coupling $\P$. 
In these ways, 
one can obtain various estimates for heat kernel, 
harmonic functions(or harmonic maps), eigenvalues etc.~by means of 
the geometry of $X$. 
These results indicate that the existence of a good coupling 
reflects the nature of $Z$ or $X$. 

Our interest in this paper is the problem of the uniqueness. 
More precisely, 
we would like to know what properties of $Z$ or $X$ are related 
to the uniqueness of couplings 
which minimize the coupling probability. 
At this moment, however, the existence of such a good coupling 
is not obvious at all in general. 
Thus we confine ourselves in a special situation 
where the existence is ensured. 

In the preceding work by E.~P.~Hsu and K.-Th.~Sturm~\cite{Hsu-Stu}, 
they discussed the uniqueness of maximal coupling 
when $X= \R^d$ and $Z$ is the Brownian motion on it. 
Motivated by the coupling inequality, 
they defined a maximal coupling 
as it minimizes the coupling probability. 
In their framework, there is a natural maximal coupling 
$\P_M \in \sC ( \P_{x_1} , \P_{x_2})$ 
called ``mirror coupling'' defined by using the reflection 
with respect to the hyperplane which maps $x_1$ to $x_2$. 
They showed that the mirror coupling 
is the unique maximal coupling 
in the class of Markovian couplings 
$\sC_0(\P_{x_1}, \P_{x_2} )$~(see \Defn{Markov}). 
They also showed by examples that 
the uniqueness no longer holds 
when we are allowed to take non-Markovian couplings. 
Their argument to derive the uniqueness uses 
the explicit form of the transition density of 
the Brownian motion. 
In this sense, their argument depends on 
the nature of the Euclidean Brownian motion. 

In order to investigate how such a uniqueness depends on the 
nature of $Z$ or $X$, 
we discuss the same uniqueness problem in 
a similar, but more general, situation. 
That is, we assume a sort of `reflection structure' 
like a reflection in Euclidean spaces 
for given initial points $x_1, x_2 \in X$. 
Then we can naturally define a mirror coupling 
$\P_M \in \sC( \P_{x_1} , \P_{x_2})$ as a maximal Markovian coupling. 
In this situation, we consider the uniqueness of 
maximal couplings in $\sC_0 ( \P_{x_1}, \P_{x_2} )$.  
As a result, 
the Brownian motion on a Riemannian manifold enjoys 
the uniqueness. 
But, as we will see, the uniqueness no longer holds 
if we consider more singular cases.
These observations show that the uniqueness is related 
to the nature of $Z$ or $X$ 
even when the mirror coupling exists. 

The organization of this paper is as follows. 
In the next section, we introduce our framework including 
the notion of `reflection structure', maximal coupling and 
Markovian coupling. 
Our main theorem gives a sufficient condition to 
the uniqueness of maximal couplings 
in $\sC_0 (\P_{x_1}, \P_{x_2} )$~(\Thm{Max}). 
We will prove \Thm{Max} in section 
\ref{sec:main} following the idea of \cite{Hsu-Stu}. 
Section \ref{sec:example} is devoted to some examples. 
On one hand, we will show that the uniqueness holds 
under the assumption on the short time asymptotic behavior of $Z$ 
and the geometry of $X$(\Thm{mfd}). 
A typical example satisfying these conditions 
is the Brownian motion on a complete Riemannian manifold 
(\Cor{Riem}). 
This framework includes the Euclidean Brownian motion 
as discussed in \cite{Hsu-Stu}. 
There we exhibit complete Riemannian manifolds 
which have the reflection structure 
with respect to specified initial points. 
On the other hand, we also show two easy examples where 
the uniqueness of maximal Markovian coupling fails 
(see~\Exam{eight} and \Exam{tree}). 
At the end of this section, we consider the case for 
the Brownian motion on $2$-dimensional Sierpinski gasket. 
We show that the uniqueness holds 
while this case is not included in the framework of \Thm{mfd}. 
In section \ref{sec:KC_couple}, 
we show that the Kendall-Cranston coupling 
coincides with our mirror coupling 
under the existence of the reflection structure. 
The Kendall-Cranston coupling is originally introduced by 
Kendall \cite{Kend} and Cranston \cite{Crans} 
for the Brownian motion 
on an arbitrary complete Riemannian manifold. 
Their coupling is useful to estimate analytic quantities 
by means of the geometric quantity such as Ricci curvature. 
But, in general, 
there is no reason why the Kendall-Cranston coupling should be  
maximal. 
The construction of their coupling is based on 
a sort of reflection of infinitesimal motion 
by means of the Riemannian geometry. 
Thus, our result is quite natural. 
It should be remarked that 
our result implies that the Kendall-Cranston coupling is 
the unique maximal coupling 
if there is a reflection structure. 

\section{Coupling of diffusions and its properties}
\label{sec:max}

Throughout this paper, we assume 
$X$ to be an arcwise-connected Hausdorff topological space 
with the second countability axiom. 
For a coupled diffusion process $(Z_1(t), Z_2(t))$, 
the coupling time $T(Z_1, Z_2)$ is defined by 
\begin{equation} \label{eq:ctime}
T(Z_1, Z_2) 
:= 
\inf 
\{
 t>0 \; ; \; Z_1 (s) = Z_2 (s) \mbox{ for all } s \ge t
\} .
\end{equation}
We set 
\begin{equation}\label{eq:cost}
\ph_t (x, y) := \frac12 \varm{x}{y}{t} .
\end{equation}
Here $\vnorm{\cdot}$ stands for the total variation norm. 
By using this function, 
the coupling inequality is written as follows: 
for every $x,y \in X$ and $\P \in \sC (\P_{x} , \P_{y} )$, 
\begin{equation} \label{eq:cineq}
 \P [ T( Z_1 , Z_2 ) > t ] \ge \ph_t(x , y) .
\end{equation}
For the proof of \eq{cineq}, it suffices to remark that 
\begin{align*}
\P [ T( Z_1 , Z_2 ) > t ] 
& \ge  
\P [ Z_1 (t) \neq Z_2 (t) ] 
\\
& \ge 
\P [ Z_1 (t) \in A , Z_2 (t) \notin A ]
\\
& \ge 
\P [ Z_1 (t) \in A ] - \P [ Z_2 (t) \in A ] 
\end{align*}
holds for arbitrary $A \in \sB(X)$. 

\bdefnn{maximal}{cf. \cite{Griff, Hsu-Stu}}
For $t>0$, we say $\P \in \sC( \P_x , \P_y )$ maximal at $t$
when the equality holds in \eq{cineq}. 
We say $\P \in \sC (\P_x , \P_y )$ maximal 
when the equality holds in \eq{cineq} for each $t>0$. 
\edefn
Let us fix $x_1 , x_2 \in X$. 
The reflection structure with respect to $x_1$ and $x_2$ 
stated in section \ref{sec:intro} means 
the following two properties assigned on $X$ and $Z$: 
\begin{description} 
\item{(A1)} \label{refl}
There is a continuous map $R : X \to X$
with $R \circ R = \mathrm{id}$ 
so that $\P_{x_1}\circ R^{-1} = \P_{x_2}$, 
\item{(A2)} \label{decomp} 
The set of fixed points $H := \{ x \in X \; ; \; R(x) = x \}$ 
separates $X$ into 
two disjoint open sets 
$X_1$ and $X_2$~(i.e., $X\setminus H = X_1 \sqcup X_2$) 
with $R(X_1)=X_2$. 
\end{description}
As an easy but significant consequence of (A2), 
every continuous path in $X$ joining $x \in X_1$ and 
$y \in X_2$ must intersect $H$. 
In general, it highly depends 
on the choice of $x_1, x_2 \in X$ 
whether (A1) and (A2) hold or not (see \Exam{torus}). 
But, we can easily verify that (A1) and (A2) are satisfied 
for the Euclidean Brownian motion for any $x_1 , x_2 \in X$. 
In that case, $R$ is an reflection with respect to a hyperplane $H$. 
Under (A1) and (A2), we can construct a mirror coupling 
of $\P_{x_1}$ and $\P_{x_2}$. 
Let $\tau := \inf \{ t >0 \; ; \; Z_1(t) \in H \}$ be 
a hitting time to $H$. 
We define the mirror coupling $\P_M$ as the law of 
$(Z_1, Z_2)$ where $Z_1$ is a copy of $(Z , \P_{x_1})$ and 
\begin{equation}
  Z_2(t) = 
  \begin{cases}
    R Z_1(t) & \mbox{if $t < \tau$} ,
    \\ 
    Z_1(t) & \mbox{if $t \ge \tau$} .
  \end{cases}
\end{equation}
By definition, 
$\P_M \in \sC (\P_{x_1} , \P_{x_2})$ and 
$\tau = T(Z_1, Z_2)$ under $\P_M$. 

\bprop{MM} 
$\P_M$ is maximal.
\eprop
For the proof, we use the following lemma. 
\blem{confined}
Suppose $\P_x \circ R^{-1} = \P_{Rx}$ for $x \in X_1$. 
Then, for each $t >0$, 
\begin{align*}
\ph_t(x, Rx) 
& = 
\P_{x} [ Z(t) \in X_1 ] - \P_{Rx} [ Z(t) \in X_1] .
\end{align*}
\elem
\bpf
By \eq{cost},
\begin{equation}   \label{eq:msup}
  \ph_t( x , Rx ) 
   = 
  \sup_{A \in \sB (X)} 
  \abra{ 
   \P_x [ Z(t) \in A ] - \P_{Rx} [ Z(t) \in A]
  } .
\end{equation}
Note that 
\begin{align*}
  \P_x [ Z(t) \in A ] - \P_{Rx} [ Z (t) \in A ]
  & = \P_x [ Z(t)\in A , \, \tau \le t ]
  - \P_{Rx} [ Z(t)\in A , \, \tau \le t ]
  \\
  & 
  \quad + \P_x [ Z(t)\in A , \, \tau > t ]
  - \P_{Rx} [ Z(t)\in A , \, \tau > t ].
\end{align*}
First we show 
\begin{equation} \label{eq:hop_step}  
\P_x [ Z(t)\in A , \, \tau \le t ]
  = \P_{Rx} [ Z(t)\in A , \, \tau \le t ]
\end{equation}
for each $A \in \sB(X)$.
By the strong Markov property, 
\begin{equation*}
  \P_x [ Z(t) \in A , \tau \le t ]
  =
  \E_x \cbra{ 
    \P_{Z (\tau)} [ Z(t-s) \in A ] |_{s= \tau}
    \; ; \; \tau\le t 
  }. 
\end{equation*}
By assumption, the law of $( Z(\tau), \tau)$  under $\P_x$ equals 
that under $\P_{Rx}$. 
Thus we have 
\begin{align*}
  \E_{x} \cbra{
    \P_{Z (\tau)} [ Z(t-s) \in A ] |_{s = \tau} 
    \; ; \; \tau \le t
  }
  &  = 
  \E_{Rx} \cbra{
    \P_{Z (\tau)} [ Z(t-s) \in A ] |_{s = \tau}
    \; ; \; \tau \le t
  }
  \\
  & = 
  \P_{Rx} \cbra{
    Z(t) \in A , \, \tau\le t
  }. 
\end{align*} 
Next, by (A2), we have 
\begin{multline*}
\P_x [ Z(t) \in A , \, \tau > t ]
- 
\P_{Rx} [ Z(t) \in A , \, \tau > t ]
\\
= 
\P_x [ Z(t) \in X_1 \cap A , \, \tau > t ]
- 
\P_{Rx} [ Z(t) \in X_2 \cap A , \, \tau > t ] .
\end{multline*}
These observations imply that 
the supremum in \eq{msup} is attained when $A=X_1$. 
\epf

\apf{\Prop{MM}}
By (A2), 
\[
\P_M [ T(Z_1, Z_2) > t ] = 
\P_{x_1} [ \tau > t ] = \P_{x_1} [ Z(t) \in X_1 , \: \tau > t ] 
- \P_{x_2} [ Z(t) \in X_1 , \: \tau > t ]. 
\]
By (A1), we can apply \eq{hop_step} for $x=x_1$. 
Thus we obtain 
\[
\P_{x_1} [ Z(t) \in X_1 , \: \tau > t ] 
- \P_{x_2} [ Z(t) \in X_1 , \: \tau > t ]
 = 
\P_{x_1} [ Z(t) \in X_1 ] 
- \P_{x_2} [ Z(t) \in X_1 ] .
\]
Hence \Lem{confined} yields the conclusion. 
\epf

\bdefn{Markov}
Let $Z^* = (Z_1 ,Z_2)$ 
be a coupling of diffusion process $Z$ 
starting from $(x_1, x_2)$ under $\P$. 
We define a canonical filtration $\{ \sF^*_t \}_{t \ge 0}$ 
by $\sF_s^* := \sg \{ Z^*(u)\; ; \; 0 \le u \le s \}$. 
We say that $\P$ is Markovian 
or $\P \in \sC_0 (\P_{x_1}, \P_{x_2})$ 
if, for each $s > 0$, the shifted process $\{ Z^*(t+s) \}_{t \ge 0}$
under $\P$ conditioned on $\sF^*_s$ is still 
a coupling of the diffusion process  
starting from $Z^*(s) = (Z_1(s), Z_2(s))$.
By using the shift operators $\{ \te_s \}_{s>0}$ defined by 
$( \te_s (Z^*) ) (t) = Z^*(s+t)$, 
$\P \in \sC_0 ( \P_{x_1} , \P_{x_2} )$ 
means  
$\P [ \;\cdot \; | \sF_s ] \circ \te_s^{-1} \in \sC (\P_{Z_1(s)}, \P_{Z_2(s)} )$ 
for each $s>0$. 
\edefn
Obviously, the mirror coupling $\P_M$ is Markovian. 
As noted in \cite{Hsu-Stu}, 
the condition that $Z^*$ is a Markovian coupling does \emph{not} imply 
that $Z^*$ is a Markov process in general.

To state our main theorem, 
we introduce a subclass of $\sC ( \P_{x_1}, \P_{x_2})$. 

\bdefn{LDP_Mirror}
We say $\P \in \hat{\sC} ( \P_{x_1}, \P_{x_2} )$ when, 
for each $t >0$, 
there is $\Xi^{(t)} \in \sB( X \times X)$ 
with $\Xi^{(t)} \subset X_1 \times X_2$ 
and 
$\P ( Z^*(t) \in ( \Xi^{(t)} )^c \cap X_1 \times X_2 ) = 0$ 
so that each $( x , y ) \in \Xi^{(t)}$ satisfies the following: 
if there is a decreasing sequence $\{ s_n \}_{n \in \N}$ 
of positive numbers with $\lim_{n\to\infty} s_n = 0$ so that 
\begin{align} \label{eq:x>y}
\P_x [ Z(s_n) \in A ] 
& \ge 
\P_y [ Z(s_n) \in A ] ,
\\ \label{eq:x<y}
\P_x [ Z(s_n) \in A' ] 
& \le 
\P_y [ Z(s_n) \in A' ] 
\end{align} 
hold for all $A \subset X_1 \cup H$, $A' \subset X_2 \cup H$ 
and all $n \in \N$, then $x = Ry$. 
\edefn 
We can easily verify $\P_M \in \hat{\sC} ( \P_{x_1}, \P_{x_2} )$. 
%
\bthm{Max}
Assume (A1) and (A2) for $x_1, x_2 \in X$. 
Let 
$\P \in \hat{\sC} ( \P_{x_1} , \P_{x_2}) \cap \sC_0 ( \P_{x_1} , \P_{x_2})$. 
If there is $t_0 >0$ so that 
$\P$ is maximal at every $t \in (0,t_0)$, 
then 
the law of $Z^*(t \wg t_0)$ under $\P$ is identical to 
that under $\P_M$. 
In particular, 
if 
$\P \in \hat{\sC} ( \P_{x_1} , \P_{x_2}) \cap \sC_0 ( \P_{x_1} , \P_{x_2})$ 
is maximal, then $\P = \P_M$. 
As a result, if 
\begin{equation} \label{eq:LDP_Mirror}
\hat{\sC} (\P_{x_1}, \P_{x_2}) \supset \sC_0( \P_{x_1}, \P_{x_2}), 
\end{equation} 
then $\P_M$ is the unique maximal coupling 
in $\sC_0 (\P_{x_1}, \P_{x_2})$. 
\ethm
\brem{ass}
(i) The conditions \eq{x>y} and \eq{x<y} 
are equivalent to the fact that 
a Hahn decomposition of 
$\P_x \circ Z(s_n)^{-1} - \P_y \circ Z(s_n)^{-1}$ 
is given by $X_1$ or $X_1 \sqcup H$ 
for each $n \in \N$. 
(ii) We can directly show that the Brownian motion on a Euclidean space 
satisfy \eq{LDP_Mirror}. 
Indeed, for every $x, y \in \R^d$, 
$\tilde{X} =\{ z \in \R^d \: ; \:  | z - x | \le | z - y | \}$ gives 
a Hahn decomposition 
of $\P_x\circ Z(s)^{-1} - \P_y \circ Z(s)^{-1}$ for each $s>0$. 
This is because the transition density depends only on the distance
for fixed $t>0$. 
(iii) In the case of Euclidean Brownian motion, more strong assertion holds: 
the maximality of $\P$ only at $t>0$ implies that 
the law of $Z^*(\cdot \wg t)$ under $\P$ is identical to 
that under $\P_M$ (see~\cite{Hsu-Stu}). 
But their proof requires some properties 
derived from 
the explicit form of the transition density of 
the Euclidean Brownian motion. 
\erem

\section{\protect{Proof of \Thm{Max}}}
\label{sec:main}

To begin with, 
we remark that (A1) produces the following auxiliary lemma. 

\blem{distant_mirror}
For each $t>0$, there is $\tilde{X}^{(t)} \in \sB(X)$ with 
$\P_{x_1} [ Z(t) \in \tilde{X}^{(t)} ]=1$ 
so that 
$\P_z \circ R^{-1} = \P_{Rz}$ for $z \in \tilde{X}^{(t)}$. 
\elem
\bpf
Take $A_i \in \sB(X)$ for $i=0, \ldots , n$ and 
$0 < s_1 < s_2 < \cdots < s_n$. 
Then the Markov property implies 
\begin{multline*}
\P_{x_1} \cbra{ 
  Z(t) \in A_0 , Z(t+s_1) \in A_1 , \ldots , Z( t+ s_n ) \in A_n 
}
\\
=
\E_{x_1} \cbra{ 
  1_{A_0} (Z(t)) \cdot \P_{Z(t)} \cbra{
    Z(s_1) \in A_1 , \ldots , Z(s_n) \in A_n 
  }
} . 
\end{multline*}
By using (A1) twice, 
\begin{align*}
\P_{x_1} [ 
  Z(t) \in A_0 , & 
  Z(t+s_1) \in A_1 , \ldots , Z( t+ s_n ) \in A_n 
]
\\
& =
\P_{x_2} \cbra{ 
  Z(t) \in R^{-1} A_0, 
  Z(t+s_1) \in R^{-1} A_1 , \ldots , Z( t+ s_n ) \in R^{-1} A_n 
}
\\
& =
\E_{x_2} \cbra{ 
  1_{ R^{-1} A_0 } ( Z(t) ) \cdot 
  \P_{Z(t)} \cbra{ 
    Z(s_1) \in R^{-1} A_1 , \ldots , Z( s_n ) \in R^{-1} A_n 
  }
}
\\
& = 
\E_{x_1} \cbra{
  1_{A_0} ( Z(t) ) \cdot 
  \P_{R Z(t)} \cbra{ 
    Z(s_1) \in R^{-1} A_1 , \ldots , Z( s_n ) \in R^{-1} A_n 
  }
}. 
\end{align*}
Since $A_0$ is arbitrary, 
there is $\tilde{X}_{s_1, \ldots, s_n ; A_1, \ldots , A_n} \in \sB(X)$
with 
$\P_{x_1} [ Z(t) \in \tilde{X}_{s_1, \ldots, s_n ; A_1, \ldots , A_n}]=1$ 
so that 
\[
\P_{x} \cbra{ 
  Z(s_1) \in A_1 , \ldots , Z( s_n ) \in A_n 
}
=  
\P_{R x} \cbra{ 
  Z(s_1) \in R^{-1} A_1 , \ldots , Z( s_n ) \in R^{-1} A_n 
} 
\]
for $ x \in \tilde{X}_{s_1, \ldots, s_n ; A_1, \ldots , A_n}$. 
Since $X$ enjoys the second countability axiom, 
there is a countable family of open sets $\mathcal{U}$ in $X$ 
so that $\sg(\mathcal{U}) = \sB(X)$. 
Thus 
\[
\tilde{X}^{(t)}  
= 
\bigcap_{ n \in \N} 
\bigcap_{
  \begin{subarray}{c}
    s_i \in \mathbb{Q} \\ 1 \le i \le n 
  \end{subarray}
}
\bigcap_{
  \begin{subarray}{c}
    A_i \in \mathcal{U} \\ 1 \le i \le n
  \end{subarray}
}
  \tilde{X}_{s_1, \ldots, s_n ; A_1, \ldots , A_n}
\]
is what we desired. 
\epf
\brem{separable}
In this paper, we used the second countability axiom of $X$ 
only for the proof of \Lem{distant_mirror}. 
Thus, if $\P_x\circ R^{-1} = \P_{Rx}$ holds for all $x\in X$, 
then $X$ need not satisfy it.  
\erem

We write $\mu_1^t = \P_{x_1} \circ Z(t)^{-1}$ and 
$\mu_2^t = \P_{x_2} \circ Z(t)^{-1}$ for simplicity. 
Let us define $\mu_0^t$ by 
\begin{equation}\label{eq:nu_zero}
\mu_0^t (A) = \mu_2^t(A \cap X_1) + \mu_1^t ( A \cap X_1^c)
\end{equation} 
for each $A \in \sB(X)$. 
By \Lem{confined}, we have $\mu_0^t \le \mu_1^t$ and $\mu_0^t \le \mu_2^t$. 
\bdefn{Trans_mirror} 
For $t>0$, 
the mirror coupling $\mu^t_M \in \sC (\mu_1^t, \mu_2^t)$ is 
the probability measure on $X \times X$ defined by 
\begin{equation}\label{eq:tmirror}
\mu_M^t ( dx dy) = \dl_{x} (dy) \mu_0^t (dx)
+ \dl_{Rx}(dy) (\mu_1^t - \mu^t_0) (dx). 
\end{equation}
\edefn
We can easily verify $\mu_M^t \in \sC ( \mu_1^t, \mu_2^t)$. 
\blem{Wasser}
Let $s, t >0$. Then for $x, y \in X$, 
\begin{equation} \label{eq:wasser_ineq}
  \inf \bbra{
    \int_{X \times X} \ph_s ( z_1 , z_2 ) \nu (dz_1 dz_2)
    \: ; \: 
    \nu \in \sC( \P_x \circ Z(t)^{-1} , \P_y \circ Z(t)^{-1} )
  }
  \ge
  \ph_{s+t} ( x , y ) .
\end{equation}
In particular, the equality holds when $(x,y)=(x_1, x_2)$. 
In this case, the infimum is attained at $\mu^t_M$. 
\elem
\bpf
Let $u_{t,E}(z) := \P_{z} [ Z(t) \in E ]$ for $E \in \sB(X)$. 
Let $\mu^t \in \sC (\P_{x}\circ Z(t)^{-1} , \P_{y}\circ Z(t)^{-1} )$.
Then 
\begin{align*}
u_{s+t,E} ( x ) - u_{s+t,E} ( y )
& = 
\E_x [ u_{s,E} (Z(t)) ] - \E_y [ u_{s,E} (Z(t)) ]
\\
& = 
\int_{X\times X}\bbra{ 
  u_{s,E} ( z_1 ) - u_{s,E} ( z_2 ) 
}
d \mu^t ( dz_1 dz_2 )
\\
& \le 
\int_{X \times X} \ph_s ( z_1 , z_2 ) d\mu^t ( dz_1 dz_2 ).
\end{align*} 
By taking the supremum on $E \in \sB(X)$ 
in the left hand side of the above inequality, 
we obtain \eq{wasser_ineq}. 
We now turn to the latter assertion. 
We set $x=x_1$ and $y= x_2$. 
By \eq{tmirror}, we have 
\begin{equation}\label{eq:meq}
\int_{X \times X}
 \ph_s( z_1, z_2 ) d\mu_M^t ( dz_1 dz_2 )
 = 
 \int_X \ph_s (z, Rz)\mu_1^t(dz) 
 - \int_X \ph_s (z, Rz) \mu_0^t (dz) .
\end{equation} 
Set $u_t(z) = u_{t,X_1}(z)$. 
Let $\tilde{X}^{(t)}$ be as in \Lem{distant_mirror}. 
By \Lem{confined}, we obtain 
\begin{equation} \label{eq:balanced}
\ph_s(z , Rz) = u_s (z) - u_s (Rz)
\end{equation}
for $z \in X_1 \cap ( \tilde{X}^{(t)} \cup R \tilde{X}^{(t)} )$. 
Thus 
\begin{equation} \label{eq:meq1}
  \int_X \ph_s (x, Rx) \mu_0^t (dx) 
   = 
  \int_{X_1} \ph_s (x, Rx) \mu_2^t(dx) + \int_{X_2} \ph_s (x, Rx) \mu_1^t(dx).
\end{equation}
Substituting \eq{meq1} to \eq{meq}, we obtain 
\begin{align*}
\int_{X \times X} \ph_s(x , y) \mu_M^t(dx dy)
& = 
\int_{X_1} \ph_s( x , Rx ) \mu_1^t(dx) 
-
\int_{X_1} \ph_s( x , Rx ) \mu_2^t(dx) 
\\
& =
\int_{X_1 \cap \tilde{X}^{(t)}} \ph_s( x , Rx ) \mu_1^t(dx) 
-
\int_{X_1 \cap R \tilde{X}^{(t)} } \ph_s( x , Rx ) \mu_2^t(dx) 
\\
& =
\int_{X_1} \bbra{ u_s( x ) - u_s( Rx ) } \mu_1^t(dx) 
-
\int_{X_1} \bbra{ u_s( x ) - u_s( Rx ) } \mu_2^t(dx) 
\\
& = \int_X u_s (x) \mu_1^t(dx) - \int_X u_s (x) \mu_2^t(dx)
\\
& = u_{s+t} (x_1) - u_{s+t}(x_2)
\\
& = \ph_{s+t} (x_1 , x_2).
\end{align*}
Here the third equality follows from \eq{balanced}.
\epf
In the following, we show 
a kind of converse assertion. 
\bprop{trans_unique} 
Let $\P \in \hat{\sC} ( \P_{x_1}, \P_{x_2} )$ and $t>0$. 
Suppose that there is a sequence $\{ s_n \}_{n \in \N}$ so that 
\begin{equation}\label{eq:wasser_eq}
\E [ \ph_{s_n} ( Z_1 (t) , Z_2 (t) ) ] = \ph_{s_n +t}(x_1 , x_2 )
\end{equation} 
holds for all $n \in \N$. 
Then $\P \circ (Z_1(t), Z_2(t))^{-1} = \mu_M^t$. 
\eprop
Let $D := \{ (x,x) \in X \times X \; ; \; x \in X \}$ 
and $\iota \: : \: X \to D$ a canonical injection. 
For the proof of \Prop{trans_unique}, we show the following lemma. 
\blem{step1}
Let $\P \in \sC ( \P_{x_1}, \P_{x_2} )$ and $t>0$. 
Suppose that there is a sequence $\{ s_n \}_{n \in \N}$ so that 
\[
\E [ \ph_{s_n} ( Z_1 (t) , Z_2 (t) ) ] = \ph_{s_n +t}(x_1 , x_2 )
\]
holds for all $n \in \N$. 
Then $\P \circ (Z_1(t), Z_2(t))^{-1} |_D = \mu_0^t \circ \iota^{-1}$. 
\elem
\bpf
Set $\mu := \P \circ ( Z_1(t) , Z_2(t) )^{-1}$. 
For simplicity, we write $\mu_i^t =: \mu_i$ for $i=0,1,2$. 
By a usual argument, 
$\mu$ is expressed in the following forms: 
\begin{equation*}
\mu(dx dy) 
 =
k_1 (x, dy ) \mu_1 ( dx )
 =
k_2 (y, dx ) \mu_2 ( dy ).
\end{equation*}
We define a coupling 
$\nu \in \sC ( \mu_1 , \mu_2 )$ by 
\begin{align*}
\nu ( dx dy ) 
& = 
\frac12 \dl_x ( dy ) \mu_0 ( dx ) 
 + 
\frac12 \int_X k_2( z, dx ) k_1 ( z, dy ) \mu_0 (dz)
\\
& \quad 
 -  
\frac12 k_1 ( x, dy ) \mu_0 ( dx )
 - 
\frac12 k_2 ( y, dx ) \mu_0 ( dy ).
 + 
\mu( dx dy )
\end{align*}
By \eq{wasser_ineq} and \eq{wasser_eq}, 
for $ s \in \{ s_n \}_{n\in \N}$, we have 
\begin{align}\nn 
0 
& \le 
\int_{X \times X} \ph_s ( x , y ) \nu ( dx dy )
-
\int_{X \times X} \ph_s ( x , y ) \mu ( dx dy )
\\ \nn
& = 
\frac12 \int_{X \times X \times X} 
\ph_s ( x , y ) k_2( z, dx ) k_1 ( z , dy ) \mu_0 ( dz )
\\ \nn
& \quad 
 - 
\frac12 \int_X \ph_s ( x , y ) k_1 ( x, dy ) \mu_0 ( dx )
 -
\frac12 \int_X \ph_s ( x , y ) k_2 ( y, dx ) \mu_0 ( dy ) 
\\ \label{eq:detour}
& = 
\frac12 \int_{X \times X \times X } 
\bbra{ \ph_s ( x , y ) - \ph_s ( z , y ) - \ph_s ( z , x ) }
k_2( z, dx ) k_1 ( z , dy ) \mu_0 ( dz ) .
\end{align}
By the triangular inequality for $\vnorm{\cdot}$, 
we have 
$
\ph_s ( x , y ) \le \ph_s ( x , z ) + \ph_s ( z , y ).
$
Thus the left hand side of \eq{detour} must be 0.
Moreover, 
there is $\Om_s \subset X \times X\times X$ with 
$\int_{\Om_s^c} k_2 ( z , dx ) k_1 ( z , dy ) \mu_0 ( dz ) =0$
so that 
\begin{equation*} 
\ph_s ( x , y ) = \ph_s ( x , z ) + \ph_s ( z , y )
\end{equation*} holds 
for each $(x, y, z) \in \Om_s$.
Note that this equality is equivalent to 
the existence of a Borel subset $E_s^{(x,y,z)} \subset X$ 
which satisfies 
\begin{equation*}
\P_x [ Z(t) \in A ] \le \P_z [ Z(t) \in A ] \le \P_y [ Z(t) \in A ]
\end{equation*}
for each Borel set $A \subset E_s^{(x,y,z)}$ 
and 
\begin{equation*}
\P_y [ Z(t) \in A' ] \le \P_z [ Z(t) \in A' ] \le \P_x [ Z(t) \in A' ]. 
\end{equation*} 
for each Borel set $A'\subset ( E_s^{(x,y,z)} )^c$. 
This fact follows from 
a simple calculation of the total variation norm 
by using Hahn decompositions. 
Let $\Om := \bigcap_{n\in \N} \Om_{s_n}$. 
We set 
\begin{align*}
A_1 
& = 
\{ ( x , y , z ) \in X \times X \times X \; ; \; x = z \},
\\
A_2 
& = 
\{ ( x , y , z ) \in X \times X \times X \; ; \; y = z \}. 
\end{align*} 
Then we claim 
\begin{equation} \label{eq:coincidence} 
\Om \subset A_1 \cup A_2 . 
\end{equation} 
Let $( z_1 , z_2 , z_3 ) \in \Om$. 
Suppose $( z_1 , z_2 , z_3 ) \notin A_1 \cup A_2$. 
Take open neighborhoods 
$V_i$ of $z_i$ ($i = 1, 2, 3$) 
with $V_i\cap V_3 = \emptyset$ for $i =1,2$. 
We choose $n \in \N$ sufficiently large so that 
$\P_{z_i} [ Z (s_n) \in V_i ] \ge 3/4$ for $i=1,2,3$. 
But, for $E_{s_n} = E_{s_n}^{(z_1,z_2,z_3)}$, we have 
\begin{align*}
\frac34 
\le 
\P_{z_3} [Z (s_n) \in V_3 ] 
& = 
\P_{z_3} [ Z (s_n) \in V_3 \cap E_{s_n} ] 
 + 
\P_{z_3} [ Z (s_n) \in V_3 \cap E_{s_n}^c ]
\\
& \le
\P_{z_2} [ Z(s_n) \in V_3 \cap E_{s_n} ]
 + 
\P_{z_1} [ Z(s_n) \in V_3 \cap E_{s_n}^c ] 
\\
& \le \frac12 .
\end{align*}
Of course it is absurd. 
Now \eq{coincidence} yields 
\begin{align*} 
\mu_0(X) 
& = 
\int_\Om k_1 ( z , dy ) k_2 ( z , dx ) \mu_0 ( dz )
\\
& \le
\int_{A_1} k_1 ( z , dy ) k_2 ( z , dx ) \mu_0 ( dz )
+ 
\int_{A_2} k_1 ( z , dy ) k_2 ( z , dx ) \mu_0 ( dz )
\\
& \quad -  
\int_{A_1 \cap A_2} k_1 ( z , dy ) k_2 ( z , dx ) \mu_0 ( dz )
\\
& = 
\mu_0 (X) 
-
\int_X 
\abra{ 
  1 -  k_1 ( z , \{ z \} ) 
} 
\abra{ 
  1 - k_2 ( z , \{ z \} ) 
}
\mu_0 ( dz ) .
\end{align*}
This equality asserts 
that 
there is $\tilde{\Om} \in \sB(X)$ with $\mu_0(\tilde{\Om}^c ) =0$ 
so that $k_1(x,\{ x \} ) =1$ or $k_2(x, \{ x \} ) =1 $ holds 
for all $x\in \tilde{\Om}$. 
Set $\tilde{\Om}_1 := \{ x \in X \; ; \; k_1 ( x , \{ x \} ) =1\}$. 
Let $\iota \: : \: X \to X \times X$ be 
given by $\iota (x) = ( x , x )$. 
For $A \in \sB(X)$, \eq{nu_zero} yields 
\begin{align*}
  \mu_0 (A) 
  & = 
  \mu_2 ( A \cap X_1 ) + \mu_1 ( A \cap X_1^c )
  \\
  & \ge
  \int_{ A \cap X_1 } k_2 ( z , \{ z \} ) \mu_2 ( dz )
  + 
  \int_{ A \cap X_1^c } k_1 ( z , \{ z \} ) \mu_1 ( dz )
  \\
  & = 
  \mu ( \iota ( A \cap X_1 ) ) 
   + 
  \mu ( \iota ( A \cap X_1^c ) )
  \\
  & = 
  \mu ( \iota ( A ) ).
\end{align*}
This estimate implies $\mu ( \iota ( \tilde{\Om}^c ) ) = 0$. 
Thus we have 
\begin{align*}
\mu ( \iota (A) ) 
& = 
\mu (\iota ( A \cap \tilde{\Om}_1 ) )
 + 
\mu (\iota ( A \cap \tilde{\Om}_1^c \cap \tilde{\Om} ) )
\\
& = 
\int_{A \cap \tilde{\Om}_1} k_1 ( z , \{ z \} ) \mu_1 ( dz )
 + 
\int_{A \cap \tilde{\Om}_1^c \cap \tilde{\Om} } 
k_2 ( z , \{ z \} ) \mu_2 ( dz )
\\
& = 
\mu_1 ( A \cap \tilde{\Om}_1 ) 
 + 
\mu_2 ( A \cap \tilde{\Om}_1^c  \cap \tilde{\Om} )
\\
& \ge  
\mu_0 ( A )
.
\end{align*}
Thus we obtain $\mu |_{D} = \mu_0 \circ \iota^{-1}$. 
\epf 
\apf{\Prop{trans_unique}} 
We use the same notation as in the proof of \Lem{step1}. 
We denote $\tilde{\mu} = \mu - \mu_0 \circ \iota^{-1}$. 
Note that $\tilde{\mu}$ is positive 
and it is absolutely continuous with respect to $\mu$ 
by \Lem{step1}. 
In order to derive $\mu = \mu^t_M$, 
we consider the integration of $\ph_s$ by $\mu^t_M$ 
for $s \in \{ s_n \}_{n \in \N}$: 
\begin{align*} 
\int_{X \times X} \ph_s ( x, y ) \mu^t_M ( dx dy )
& =
\int_X \ph_s ( x, Rx ) \mu_1 (dx) - \int_X \ph_s ( x, Rx ) \mu_0 (dx)
\\
& =
\int_{X \times X} \ph_s ( x, Rx ) \mu ( dx dy) 
- \int_{X \times X} \ph_s ( x, Rx ) \mu_0 \circ \iota^{-1}(dx dy)
\\
& = 
\int_{X \times X} \ph_s ( x, Rx ) \tilde{\mu} ( dx dy ).
\end{align*}
By virtue of (A1), we also obtain 
\[
\int_{X \times X} \ph_s ( x, y ) \mu^t_M ( dx dy )
 = 
\int_{X \times X} \ph_s ( y, Ry ) \tilde{\mu} ( dx dy ).
\]
By \Lem{Wasser}, 
\begin{align} \nn
  0 
  & =
  \int_{X \times X} \ph_s ( x , y ) \mu^t_M ( dx dy )
  -
  \int_{X \times X} \ph_s ( x , y ) \mu ( dx dy )
  \\ \label{eq:compare_mirr}
  & = 
  \frac12 
  \int_{ X \times X } 
  \bbra{
    \ph_s ( x , Rx )
    + 
    \ph_s ( y , Ry )
    -
    2 \ph_s ( x , y ) 
  }
  \tilde{\mu} ( dx dy ). 
\end{align}
By (A1) and \eq{nu_zero}, 
$\tilde{\mu} ( X_1^c \times X ) = \tilde{\mu} ( X \times X_1 ) = 0$. 
This fact together with \eq{compare_mirr} yields 
\begin{equation}
\int_{ X_1 \times X_1^c } 
\bbra{
  \ph_s ( x , Rx )
  + 
  \ph_s ( y , Ry )
  -
  2 \ph_s ( x , y ) 
}
\tilde{\mu} ( dx dy ) = 0. 
\end{equation}
Let $\tilde{X}^{(t)}$ be as given in \Lem{distant_mirror}. 
For $z \in X_1 \cap \tilde{X}^{(t)} $, 
\begin{align*}
  \ph_s ( z , Rz ) 
  &  = 
  \P_{z} [ Z(s) \in X_1 ] - \P_{Rz} [ Z(s) \in X_1 ] 
  \\
  & =
  \P_{z} [ Z(s) \in X_1 ] - \P_{z} [ Z(s) \in X_2 ] .
\end{align*}
Thus, 
for 
$
x, y \in \abra{ X_1 \times X_1^c } 
 \cap 
( \tilde{X}^{(t)} \times R \tilde{X}^{(t)} )
$, 
\begin{align} \nn
  \ph_s ( x , Rx )
  & + 
  \ph_s ( y , Ry )
  -
  2 \ph_s ( x , y )  
  \\ \nn
  & =
  \P_{x} [ Z(s) \in X_1 ] - \P_{y} [ Z(s) \in X_1 ] 
  +
  \P_{y} [ Z(s) \in X_2 ] - \P_{x} [ Z(s) \in X_2 ] 
  \\ \nn
  & \quad 
  - 2 \sup_{A \in \sB(X)} | \P_x[ Z(s) \in A ] - \P_y [Z(s) \in A] |  
  \\ \label{eq:compar_mirr2}
  & \le 0 .
\end{align} 
Note that 
\[
\mu ( ( \tilde{X}^{(t)} \times R \tilde{X}^{(t)} )^c )
 \le
\mu_1 ( ( \tilde{X}^{(t)} )^c )
 + 
\mu_2 ( ( R \tilde{X}^{(t)} )^c )
 = 0 
\]
since $\mu \in \sC ( \mu_1 ,\mu_2 )$. 
Hence there is $\tilde{E}_s \subset \abra{ X_1 \times X_1^c } 
\cap ( \tilde{X}^{(t)} \times R \tilde{X}^{(t)} )$ 
with $\tilde{\mu} (\tilde{E}^c_s )=0$ so that 
the equality holds in \eq{compar_mirr2} 
for $(x, y ) \in \tilde{E}_s$. 
Let $E = \Xi^{(t)} \cap ( \bigcap_{n\in \N} \tilde{E}_{s_n} )$. 
Here $\Xi^{(t)}$ is given in \Defn{LDP_Mirror} associated with 
$\P \in \hat{\sC} ( \P_{x_1} , \P_{x_2} )$. 
Then $\tilde{\mu} ( E^c )  = 0$
and 
\begin{align*}
  \frac12 \varm{x}{y}{s_n}
  & = 
  \P_{x} [ Z(s_n) \in X_1 ] - \P_{y} [ Z(s_n) \in X_1 ] 
  \\ 
  & =
  \P_{y} [ Z(s_n) \in X_2 ] - \P_{x} [ Z(s_n) \in X_2 ] 
\end{align*} 
for all $( x, y ) \in E$ and $n \in \N$. 
Hence the property of $\Xi^{(t)}$ 
immediately implies $x = Ry$ for every $(x,y) \in E$ 
(cf. \Rem{ass}(i)). 
It yields $\mu = \mu^t_M$. 
\epf

\apf{\Thm{Max}}
Let $t_0 \in (0, \infty]$ and  
$\P \in \hat{\sC} ( \P_{x_1} \P_{x_2} ) \cap \sC_0 ( \P_{x_1}, \P_{x_2})$ maximal 
at each $t \in (0,t_0)$. 
Note that 
\begin{equation} \label{eq:hit=opt}
T(Z_1, Z_2) 
 = 
\inf \{ 
s > 0 \; ; \; Z_1 (s) = Z_2 (s) 
\}
\quad \mbox{$\P$-a.s.}
\end{equation}
holds since $\P$ is maximal. 
Take $s,t >0$ with $s+t < t_0$. 
By the maximality of $\P$ at $s+t$, 
\begin{align*}
\ph_{s+t} ( x_1 , x_2 ) 
& =
\P [ T ( Z_1 , Z_2 ) > s+t ]
\\
& =
\E [\, \P [ T ( Z_1  , Z_2 ) > s + t \, | \, \sF^*_t ] \, ]  .
\end{align*}
Since $\P \in \sC_0 ( \P_{x_1} , \P_{x_2} )$, 
\eq{cineq} yields 
$
  \P [ \, T ( Z_1 , Z_2 ) > s + t \, | \, \sF^*_t \, ]
  \ge \ph_s ( Z_1(t) , Z_2(t) ).
$
In addition, by \eq{wasser_ineq}, 
$
  \E [ \ph_s ( Z_1 (t) , Z_2 (t) ) ]
  \ge 
  \ph_{s+t}( x_1 , x_2 ).
$
Hence we obtain 
\begin{equation*}
  \E [ \ph_s ( Z_1 (t) , Z_2 (t) ) ] = \ph_{s+t} ( x_1 , x_2 ).
\end{equation*}
Letting $s \to 0$, \Prop{trans_unique} yields 
$\P \circ ( Z_1 (t) , Z_2 (t) )^{-1} = \mu_M^t$. 
Since $t \in  (0, t_0)$ is arbitrary, 
it implies that 
\[
\P [ Z_2 (t) = Z_1 (t) \mbox{ or } Z_2 (t) = R Z_1 (t) 
\mbox{ for all } t \in (0,t_0) ] =1.
\]
Recall that $\tau$ is the first hitting time of $Z_1$ to $H$. 
The above equality implies that $\tau$ equals 
the first hitting time of $Z_2$ to $H$ $\P$-almost surely. 
In addition, by (A2), for each $t \in (0,t_0)$, 
\begin{equation} \label{eq:before_couple}
\{ t \le \tau \} \subset \{ Z_2 (t) = R Z_1 (t) \} 
\quad \P\mbox{-a.s.}.
\end{equation}
Thus it suffices to show that 
\begin{equation} \label{eq:coupled}
\{ \tau < t \} \subset \{ Z_2 (t) = Z_1 (t) \}
\quad \P\mbox{-a.s.}.
\end{equation}
Note that \eq{tmirror} implies 
$
\P [ Z_1 (t) = Z_2 (t) ] = \mu_0^t (X) .
$
By the maximality of $\P$, \Lem{confined} and \eq{nu_zero}, 
\begin{align*}
\P [ T(Z_1, Z_2) \le t ] = 1 - \ph_t ( x_1 , x_2 ) 
& = 
1 - \P_{x_1} [ Z(t) \in X_1 ] + \P_{x_2} [ Z(t) \in X_1 ]
\\
& = \P_{x_1} [ Z(t) \in X_1^c ] + \P_{x_2} [ Z(t) \in X_1 ]
\\
& = 
\mu_0^t ( X ) .
\end{align*}
It means 
\begin{equation} \label{eq:prob_coupled}
\P [ Z_1 (t) = Z_2 (t) ] = \P [ T(Z_1, Z_2) \le t ].
\end{equation} 
Since we have \eq{before_couple}, 
\begin{align*}
  \{ Z_1 (t) = Z_2 (t) \} & \subset \{ \tau  < t \} 
& \hspace{-3cm} \P\mbox{-a.s.},
\\
\{ \tau < t \} & \subset \{ T(Z_1 , Z_2) < t \}
& \hspace{-3cm} \P\mbox{-a.s.}. 
\end{align*}
The second inclusion follows from \eq{hit=opt}. 
Combining them with \eq{prob_coupled}, we obtain \eq{coupled} and 
it completes the proof.
\epf

\section{Examples and counterexamples}
\label{sec:example}

Let us consider several examples of $X$ and $Z$ with 
(A1) and (A2) for given $x_1,x_2 \in X$.
First we state a sufficient condition 
for \eq{LDP_Mirror} to be satisfied. 
A key ingredient is 
the Varadhan type short time asymptotic behavior 
of transition probabilities \eq{Varadhan}. 
In order to state it in a general form, 
we introduce some terms concerning the metric geometry. 
Let $(X,d)$ be a metric space. 
We call a curve $\gm \: : \: [0,1] \to X$ geodesic 
if, for each $t \in [0,1]$, there exists $\dl >0$ so that
$d ( \gm (t) ,\gm (s) ) = | t - s | d( \gm (0) , \gm (1) )$ 
for $|t - s| < \dl$. 
Recall that a metric space $(X,d)$ is geodesic when,  
for each $x, y \in X$, there is a rectifiable curve in $X$ 
whose length realizes $d(x,y)$. 
Note that such a curve always becomes a geodesic 
by a suitable re-parameterization. 
We call it a minimal geodesic joining $\gm (0)$ and $\gm (1)$. 
A geodesic metric space $(X,d)$ is said to be non-branching 
when, for any two geodesics $\gm$ and $\gm'$ 
of the same length with $\gm (0) = \gm' (0)$, we have 
$\inf \bbra{ t > 0  ; \: \gm (t) \neq \gm' (t) } = 0 \mbox{ or }\infty$.
Here we follow the usual manner $\inf \emptyset = \infty$.
\bthm{mfd}
Let $(X,d)$ be a non-branching geodesic metric space 
and 
$( Z , \P_x )$ 
a diffusion process on it.
Suppose that 
(A1) and (A2) hold for given $x_1 , x_2 \in X$. 
In addition, we assume the following for each $t>0$: 
\begin{enumerate}
\item
the support of the law of $Z(t)$ under $\P_{x_1}$ equals $X$, 
\item
there exist 
\begin{itemize}
\item
an increasing function $\ro \: : \: (0,\infty) \to (0,\infty)$ 
with $\lim_{s\to 0} \ro(s) = 0$, 
\item
a strictly increasing function $\Psi \: : \: [0, \infty) \to [0,\infty)$, 
\item
a sequence $\{ s_n \}_{n \in \N}$ of positive numbers 
with $\lim_{n\to\infty} s_n =0$ and 
\item
$Y^{(t)} \in \sB(X)$ 
with $RY^{(t)} = Y^{(t)}$ and $\P_{x_1} [ Z(t) \notin Y^{(t)} ] =0$ 
\end{itemize} 
so that 
\begin{equation}
  \label{eq:Varadhan}
  - \lim_{n\to \infty} \ro(s_n) \log \P_x \cbra{ Z(s_n) \in A }
  = 
  \Psi ( d(x, A) )
\mbox{ for each } A \in \sB(X)  
\end{equation} 
holds for $x \in Y^{(t)}$.
\end{enumerate}
Then \eq{LDP_Mirror} holds. 
\ethm
\bpf 
Take $z \in \tilde{X}^{(t)}\cap Y^{(t)}$ for $t>0$, 
where $\tilde{X}^{(t)}$ is as in \Lem{distant_mirror}. 
Then 
$\P_z \cbra{ Z(s) \in A } = \P_{Rz}\cbra{ Z(s) \in R A}$ holds 
for each $A \in \sB(X)$. 
Thus \eq{Varadhan} yields $d(z, A) = d (Rz , RA)$.
By taking $A=B_r(w)$, 
a ball of radius $r>0$ centered at $w \in X$, 
and taking $r\to 0$, 
we obtain $d(z,w)=d(Rz,Rw)$. 
Since $R$ is continuous, the condition (i) implies that 
$R$ acts on $X$ as isometry. 
Let $x \in Y^{(t)}\cap X_1$ and $y \in Y^{(t)} \cap X_2$. 
Suppose \eq{x>y} and \eq{x<y} holds. 
Take $z \in X_1$. 
Since $X_1$ is open, $B_r(z) \subset X_1$ 
for all sufficiently small $r>0$. 
For such $r$, 
\eq{Varadhan} and \eq{x>y} implies that,  
\begin{align}
  \nn
  \Psi (  d ( x , B_r(z) ) ) 
  & =
  - \lim_{n\to \infty} \ro(s_n) \log \P_x [ Z(s_n) \in B_r(z) ]
  \\ \nn
  & \le 
  - \lim_{n\to \infty} \ro(s_n) \log \P_y [ Z(s_n) \in B_r(z) ]
  \\   \label{eq:x,z>z,y}
  & =
  \Psi ( d (y, B_r(z) ) ).
\end{align}
It immediately implies $d(x, B_r(z)) \le d (y , B_r(z))$. 
By taking $r \to 0$, we obtain $d(x, z) \le d(y, z)$.
In the same way, for $z\in X_2$, we obtain $d(x, z) \ge d(y, z)$.
These two estimates yield 
\begin{equation} \label{eq:same_distance}
d(x,z) = d(y,z) \mbox{ for } z \in H. 
\end{equation} 
Take a minimal geodesic $\gm_0 \: : \: [0,1] \to X$ 
joining $x$ and $y$. 
By the remark after (A2), 
there exists $t_0 \in [0,1]$ so that $\gm_0(t_0) \in H$. 
By \eq{same_distance}, 
we have $d(x, \gm_0(t_0) ) = d(y, \gm_0(t_0))$. 
In addition, $t_0 =1/2$ follows. 
Again by \eq{same_distance}, 
$d(x, \gm_0(1/2)) = d(x, H)$ holds. 
Let $\gm_1$ be a curve joining $x$ and $Rx$ given 
by 
\[
\gm_1 (t) = 
\begin{cases}
\ds 
\gm_0(t) 
 & 
\mbox{if }t \in \cbra{ 0 , 1/2 } ,
\\
\ds 
R ( \gm_0(1-t) ) 
 & 
\mbox{if }t \in \left( 1/2 , 1 \right] .
\end{cases}
\]
Then 
$\gm_1$ is a minimal geodesic joining $x$ and $Rx$ 
because $d(x, \gm_1(1/2) ) = d ( Rx, \gm_1(1/2) ) = d (x, H)$.
Since $X$ is non-branching, we obtain $Rx=y$. 
Thus, once we set 
$\Xi^{(t)} = ( Y^{(t)} \cap X_1 ) \times ( Y^{(t)} \cap X_2 )$,  
$\P [ Z^* (t) \in ( \Xi^{(t)} )^c \cap X_1 \times X_2 ] =0$ holds 
for each $\P \in \sC ( \P_{x_1} , \P_{x_2} )$. 
This means $\hat{\sC} ( \P_{x_1} , \P_{x_2}) = \sC ( \P_{x_1}, \P_{x_2} )$ and 
therefore \eq{LDP_Mirror} holds. 
\epf
\brem{density} 
If our diffusion process $Z(t)$ has 
a continuous transition density $p_t(x,y)$ 
with respect to a Radon measure $m$, 
that is, $\P_x[ Z(t) \in A ] = \int_A p_t (x, y) m(dy)$, 
Then \eq{Varadhan} in \Thm{mfd} is replaced as follows: 
\begin{equation} \label{eq:Varadhan2}
  - \lim_{n\to \infty} \ro(s_n) \log p_{s_n} (x, y)
  = 
  \Psi ( d(x, y) )
\mbox{ for each } y \in X.
\end{equation}
Indeed, \eq{x>y} and \eq{x<y} imply that 
$p_{s_n}(x,z) \ge p_{s_n}(y,z)$ for $z \in X_1$ and 
$p_{s_n}(x,z) \le p_{s_n}(y,z)$ for $z \in X_2$. 
Thus the same proof works. 
\erem
\bcor{Riem}
Let $X$ be a complete Riemannian manifold and 
$Z(t)$ the Brownian motion on $X$. 
Assume $X$ to satisfy (A1) and (A2). 
Then \eq{LDP_Mirror} follows. 
\ecor
\bpf
In this case, $Z$ has a continuous transition density $p_t(x,y)$.  
Letting $\ro(s) := 2s$, $\Psi(v) = v^2$
and any sequence $\{s_n\}_{n \in \N}$ with $\lim_{n\to \infty} s_n =0$, 
\eq{Varadhan2} follows 
from \cite{Norris} 
for every $x\in X$ (see~\Rem{density}). 
The condition (i) in \Thm{mfd} 
comes from strict positivity of the transition density. 
It is well-known that 
all properties imposed on $X$ in \Thm{mfd} hold.
\epf
We can also apply \Thm{mfd} to Alexandrov spaces. 
These metric spaces are an generalization of a complete Riemannian manifold 
with sectional curvature bounded below~(see \cite{Bur_Gro_Per} for details).
\bcor{Alex}
Let $(X,d)$ be an Alexandrov space and 
$Z(t)$ a diffusion process on $X$ corresponding to 
a canonical regular Dirichlet form on $X$ 
constructed in \cite{Kuw_Mac_Shi}~(see \cite{Kuw_Shi03} also). 
Assume $X$ to satisfy (A1) and (A2). 
Then \eq{LDP_Mirror} follows. 
\ecor
\bpf
In this case, there is a continuous transition density $p_t(x,y)$ of $Z$. 
Letting $\ro(s) = 2s$, $\Psi(v) = v^2$
and any sequence $\{s_n\}_{n \in \N}$ with $\lim_{n\to \infty} s_n =0$, 
\eq{Varadhan2} follows 
from Corollary 2 of \cite{Renes_compar} for every $x\in X$
(see~\Rem{density}). 
As in the case of Riemannian manifolds, the condition (i) follows 
from positivity of the transition density. 
By definition, $X$ is a geodesic space. 
The curvature condition on $X$ easily implies 
that $X$ is non-branching. 
Thus we can apply \Thm{mfd}. 
\epf
\brem{totally_geodesic}
Let $X$ be a Riemannian manifold and $Z$ the Brownian motion on it. 
(i) If (A1) and (A2) are satisfied for given initial points, 
then the argument in the proof of \Thm{mfd} implies that $R$ is isometry. 
In this case, $H$ is a totally geodesic 
smooth submanifold of $X$~(see \cite{Kob-NomII} p.61, for example). 
In particular, $H$ becomes a complete Riemannian manifold. 
In addition, $H$ is of codimension 1. 
(ii) If (A2) are satisfied with respect to an isometry $R$ 
with $R \circ R = \mbox{id}$, then (A1) follows for each $x_1 , x_2 \in X$ 
with $R x_1 = x_2$. 
\erem
In what follows, we will see some manifolds
satisfying the conditions (A1) and (A2). 
In all cases, we assume $Z$ to be the Brownian motion. 
\bexam{spaceform}
We consider the case $X$ is 
an irreducible Riemannian global symmetric space of constant curvature. 
We will review that (A1) and (A2) are satisfied 
for every distinct pair $x_1, x_2 \in X$ of starting points 
in these cases. 
By \Rem{totally_geodesic} (ii), It suffices to find 
an isometry $R$ with $R\circ R = \mbox{id}$, $R x_1 = x_2$ 
satisfying (A2). 
The flat case, i.e. $X=\R^n$, is considered in \cite{Hsu-Stu}. 

In the case of positive curvature, $X$ is a sphere: 
\begin{equation*} 
X = \mathbb{S}^n = 
\bbra{ 
  z =(z_0, z_1, \ldots, z_n) \in \R^{n+1} 
  \: ;\: 
  z_0^2 + \cdots + z_n^2  = r 
} 
\end{equation*}
with a metric induced from the canonical metric on $\R^{n+1}$.
Take $x_1, x_2 \in X$ with $x_1 \neq x_2$. 
Then we can easily verify that 
the restriction of the reflection in $\R^{n+1}$ 
with respect to a hyperplane 
fulfills all of our requirements. 

In the case of negative curvature, $X$ is a hyperbolic space: 
\[
X = \mathbb{H}^n = 
\bbra{ 
  z =(z_0, z_1, \ldots, z_n) \in \R^{n+1} 
  \: ;\: 
  - z_0^2 + z_1^2 + \cdots + z_n^2  = -r , z_0 >0
} 
\] 
with a metric induced from the Lorentz metric on $\R^{n+1}$. 
Take $x_1, x_2 \in X$ with $x_1 \neq x_2$. 
Let $m$ be the midpoint of $x_1$ and $x_2$. 
By homogeneity, we may assume $m= (r, 0, \ldots, 0)$. 
By arranging the chart appropriately, we may assume 
$x_1 = ( z_0, z_1, 0, \ldots, 0)$. 
Then $x_2 = (z_0, -z_1 , 0 , \ldots, 0)$. 
Set 
$
R\: : \: (z_0, z_1, \ldots, z_n ) \mapsto (z_0 , -z_1, z_2 , \ldots ,z_n).
$
Then $R$ fulfills all of our requirements. 
\eexam

\brem{const_curv}
The converse of \Exam{spaceform} is true 
for an irreducible Riemannian global symmetric space $X$ 
in the following sense. 
If there exists an isometry $R$ satisfying (A1) and (A2) 
for some pair $x_1, x_2 \in X$, 
then $X$ must be of a constant curvature. 
It follows from the result in \cite{Iwah} 
(cf.~\Rem{totally_geodesic}(i)). 
\erem
\begin{figure}
\begin{center}
\begin{pspicture}(-1.3,-1.525)(2.8,2.575)
\psline{->}(0,-1.025)(0,2.575)
\psline{->}(-1.3,0)(2.8,0)
\rput(0.14,0.25){$o$}
\psdot(0.7,0)
\rput(0.7,-0.3){$x_1$}
\rput(0,0.42){\rnode{C}{\psdot}}
\rput(-0.8,0.72){\rnode{D}{$x_2$}}
\ncarc{->}{D}{C}
\psline(-0.3,-0.525)(1.8,-0.525)(1.8,1.575)(-0.3,1.575)(-0.3,-0.525)
\psline{->}(0.75,-0.525)(0.76,-0.525)
\psline{->}(0.75,1.575)(0.76,1.575)
\psline{->}(-0.3,0.525)(-0.3,0.535)
\psline{->}(1.8,0.525)(1.8,0.535)
\psdot(-0.240625,-0.525)
\rput(-0.340625,-0.755){$z_1$}
\psdot(0.940625,1.05)
\rput(0.580625,1.05){$z_2$}
\psdot(-0.240625,1.575)
\rput(-0.340625,1.805){$z_3$}
\psline[linewidth=1.5pt](-0.240625,-0.525)(0.940625,1.05)(-0.240625,1.575)
\psdot(1.6953125,-0.525)
\rput(1.7953125,-0.755){$z_4$}
\psdot(1.1046875,1.05)
\rput(1.40046875,1.05){$z_5$}
\psdot(1.6953125,1.575)
\rput(1.7953125,1.805){$z_6$}
\psline[linewidth=1.5pt](1.6953125,-0.525)(1.1046875,1.05)(1.6953125,1.575)
\rput(0.75,-1.525){$a=1/3,b=1/5$}
\end{pspicture}
\hspace{1cm}
\begin{pspicture}(-1.63,-1.63)(2.47,2.47)
\psline{->}(0,-1.13)(0,2.47)
\psline{->}(-1.63,0)(2.47,0)
\rput(-0.26,0.2){$o$}
\psdot(0.42,0)
\rput(0.42,-0.3){$x_1$}
\psdot(0,0.42)
\rput(-0.3,0.62){$x_2$}
\psline(-0.63,-0.63)(1.47,-0.63)(1.47,1.47)(-0.63,1.47)(-0.63,-0.63)
\psline{->}(0.42,-0.63)(0.43,-0.63)
\psline{->}(0.42,1.47)(0.43,1.47)
\psline{->}(-0.63,0.42)(-0.63,0.43)
\psline{->}(1.47,0.42)(1.47,0.43)
\psdot(-0.63,-0.63)
\rput(-0.83,-0.83){$z_1$}
\rput(1.05,1.05){\rnode{A}{\psdot}}
\psdot(-0.63,1.47)
\rput(-0.83,1.67){$z_3$}
\psline[linewidth=1.5pt](-0.63,-0.63)(1.05,1.05)(-0.63,1.47)
\psdot(1.47,-0.63)
\rput(1.67,-0.83){$z_4$}
\psdot(1.47,1.47)
\rput(1.67,1.67){$z_6$}
\psline[linewidth=1.5pt](1.47,-0.63)(1.05,1.05)(1.47,1.47)
\rput(0.8,1.97){\rnode{B}{$z_2=z_5$}}
\ncarc{->}{B}{A}
\rput(0.42,-1.63){$a=b=1/5$}
\end{pspicture}
\hspace{1cm}
\begin{pspicture}(-1,-1)(3.1,3.1)
\psline{->}(0,-0.5)(0,3.1)
\psline{->}(-1,0)(3.1,0)
\rput(-0.2,-0.2){$o$}
\psdot(0,0)
\rput(-0.3,0.3){$x_2$}
\psdot(0.7,0)
\rput(0.7,-0.3){$x_1$}
\psline(0,0)(2.1,0)(2.1,2.1)(0,2.1)(0,0)
\psline{->}(1.05,0)(1.06,0)
\psline{->}(1.05,2.1)(1.06,2.1)
\psline{->}(0,1.05)(0,1.06)
\psline{->}(2.1,1.05)(2.1,1.06)
\psline[linewidth=1.5pt](0.35,0)(0.35,2.1)
\psline[linewidth=1.5pt](1.4,0)(1.4,2.1)
\rput(1.05,-1){$a=1/3,b=0$}
\end{pspicture}
\begin{minipage}{0.48\textwidth}
\begin{center}
\vspace{1cm}
Fig.1 
\end{center}
\end{minipage}
\end{center}
\end{figure}
\bexam{torus}
Let us consider 2-dimensional torus 
$X = \mathbb{T}^2 = ( \R / \!\!\sim )^2$, 
where $\sim$ identifies $x$ with $x+n$ 
for each $x \in \R$ and $n \in \N$. 
Let $\pi \: :\: \R^2 \to \mathbb{T}^2$ be the canonical projection. 
We denote $\pi(x)$ by $[x]$. 
Take $x_1,x_2 \in X$ with $x_1 \neq x_2$. 
By arranging an appropriate chart, 
we may assume that $x_1 =[ (a, 0) ]$ and $x_2 =[ (0,b) ]$ 
for $0 \le b \le a \le 1/2$. 
Let $K = \bbra{ z \in \mathbb{T}^2 \: ; \: d(z,x_1) = d(z, x_2) }$. 
If (A1) and (A2) are satisfied for $x_1, x_2 \in X$, 
then $R$ is isometry and 
$H \subset K$ must hold (cf. \Rem{totally_geodesic}). 
In the following, we will write $K$ explicitly. 
First we consider the case $b \neq 0$. 
Take six points $z_1, z_2, z_3, z_4, z_5, z_6 \in \R^2$ as follows:
\begin{align*}
z_1 
& = 
\abra{
  \frac{1}{2a}\abra{ 
    a^2 + b^2 -b  
  }
  ,
  b - \frac12
},
\\
z_2 
& = 
\abra{
  \frac{1}{2a}\abra{ 
    a^2 - b^2 + b 
  }  
  ,
  \frac12 
},
\\
z_3 
& = 
\abra{
  \frac{1}{2a}\abra{ 
    a^2 + b^2 -b 
  }
  ,
  b+\frac12
},
\\
z_4 
& = 
\abra{
  \frac{1}{2(1-a)}
  \abra{
    - a^2 - b^2  + b + 1 
  }
  ,
  b-\frac12
},
\\
z_5
& = 
\abra{
  \frac{1}{2(1-a)}
  \abra{ 
    - a^2 + b^2 - b + 1 
  }
  ,
  \frac12
},
\\
z_6 
& = 
\abra{
  \frac{1}{2(1-a)}
  \abra{
    - a^2 - b^2  + b + 1 
  }
  ,
  b+\frac12
} .
\end{align*}
Let $l_{ij}$ be a line segment in $\R^2$ 
whose endpoints are $z_i$ and $z_j$. 
Then 
$K = \pi ( l_{12} \cup l_{23} \cup l_{45} \cup l_{56} )$ 
holds.
We can easily verify that $K$ has singularity 
at $[z_2]$ or $[z_5]$~(see Fig.1)  
and 
$H$ cannot be contained in $K$ by \Rem{totally_geodesic}(i). 
Thus there is no reflection structure. 
Next we consider the case $b=0$. 
Then we have 
$
K =
\pi 
\abra{
  \{ (a/2, q)\: ;\: q \in [0,1] \} 
  \cup 
  \{ ((1+a)/2, q) \: ;\: q \in [0,1] \} 
}
$ (see~Fig.1). 
Thus a map $R \: : \: \mathbb{T}^2 \to \mathbb{T}^2$ defined
by $R( [(p,q)] )  = [(a-p,q)]$ satisfies (A1) and (A2) 
with $H=K$.  
\eexam 
\bexam{rot_symm} 
Let $X$ be a complete Riemannian manifold given 
by the direct product of two manifolds $Y_1$ and $Y_2$. 
We assume that the Riemannian metric on $X$ has 
a form $h(y_2) \sg_1(dy_1\otimes dy_1) + \sg_2(dy_2\otimes dy_2)$, where 
$h$ is a positive function on $Y_2$ and 
$\sg_i$ is a Riemannian metric on $Y_i$. 
We also assume (A1) and (A2) on $(Y_1, \sg_1)$ 
for given starting points $y_1^{(1)},y_2^{(1)} \in Y_1$. 
Then we can extend the reflection structure on $Y_1$ 
to $X$ in a natural way. 
As the result, for any $y^{(2)} \in Y_2$, 
(A1) and (A2) are satisfied for 
$(y_1^{(1)}, y^{(2)} ), (y_2^{(1)}, y^{(2)} ) \in X$. 

The same argument works for rotationally symmetric manifolds. 
Take a function $h \: : \: (0,\infty) \to (0,\infty)$ 
which has a smooth extension to $[0,\infty)$ and 
satisfies $h(0)=0$ and $h'(0)=1$. 
By using $h$, we define a metric on $(0,\infty)\times \mathbb{S}^{n-1}$ 
by $ds^2 = dr^2 + h(r) d\te^2$ 
for $(r,\te) \in (0,\infty)\times \mathbb{S}^{n-1}$. 
Let $X$ be the completion of $(0,\infty) \times \mathbb{S}^{n-1}$ 
by adding one point $o$. 
$o$ is the limit of $(r,\te)$ as $r\to 0$ 
for each $\te\in \mathbb{S}^{n-1}$. 
Take $x_1 =(r , \te_1)$ and $x_2 = (r, \te_2)$ for 
some $r \in (0,\infty)$ and $\te_1 , \te_2 \in \mathbb{S}^{n-1}$. 
Then, the Brownian motion with starting points $x_1$ and $x_2$ 
satisfies (A1) and (A2). 
\eexam
Next we give two simple examples satisfying (A1) and (A2) 
while the uniqueness of maximal Markovian coupling does not hold. 
\begin{figure} 
\begin{center}
\begin{pspicture}(6.2,4)
\psline[linestyle=dashed](3,3.5)(3,0.5)
\scalebox{1}[.4]{
  \psarc[linewidth=1pt]{->}(3,2){0.4}{130}{410}
}
\rput(3.6,0.5){$R$}
\psline[linestyle=dashed](2,2)(6,2)
\scalebox{.4}[1]{
  \psarcn[linewidth=1pt]{->}(14,2){0.4}{140}{220}
}
\rput(5.8,2.6){$\h$}
\pscircle(2,2){1}
\pscircle(4,2){1}
\psdot[dotsize=3pt 3](1,2) \rput(0.7,2.2){$x_1$} 
\psdot[dotsize=3pt 3](5,2) \rput(4.6,1.8){$x_2$} 
\rput(5,3.5){\rnode{A}{$z_1=z_2$}} 
\rput(3.1,2){\rnode{B}{\psdot[dotsize=3pt 3](-0.1,0)}}
\ncarc{->}{A}{B}
\end{pspicture}
\hspace{3cm}
\begin{pspicture}(-2.5,-2.5)(2.5,2.5)
\psdot(0,0)
\rput(.3,-.3){$p_0$}
\rput(1.4,0){$p_1$}
\rput(-0.7,0.6){$p_2$}
\rput(-0.7,-0.6){$p_3$}
\rput(1.9,-1.1){$p_{11}$}
\rput(1.9,1.1){$p_{12}$}
\rput(0.4,1.8){$p_{21}$}
\rput(-1.6,0.6){$p_{22}$}
\rput(-1.6,-1.1){$p_{31}$}
\rput(0.4,-1.8){$p_{32}$}
\psarc[linewidth=.8pt]{<->}(0,0){1.5}{100}{140}
\rput(-1.1,1.4){$\h$}
\psarc[linewidth=.8pt]{<->}(0,0){.8}{10}{110}
\rput(0.7,0.7){$R$}
\SpecialCoor \degrees[6]
\psline[linestyle=dashed](2;1)
\psline[linestyle=dashed](2;2)(1;2)
\multido{\i=0+2}{3}{
 \psline(1;\i)
}
\multido{\ia=1+2,\ib=3+2}{3}{
  \psset{origin={1;\ia}}
  \psdot(0,0)
  \multido{\i=\ib+2}{2}{
    \psline(1;\i)
    \psdot(1;\i)
  }
}
\end{pspicture}
\end{center}
\begin{minipage}{0.48\textwidth}
\begin{center}
Fig.2 
\end{center}
\end{minipage}
\hfill
\begin{minipage}{0.48\textwidth}
\begin{center}
Fig.3 
\end{center}
\end{minipage}
\end{figure}

\bexam{eight}{\bf (Fig.2)}
Let $Y_1$ and $Y_2$ be two copies of $\mathbb{T}^1 = [0,1]/\!\!\sim$.
For $z_i=1/2 \in Y_i$ ($i=1,2$), we set 
$X = Y_1 \sqcup Y_2 / \!\!\sim$ where $\sim$ means 
the identification of $z_1$ and $z_2$. 
Set $x_1 =0 \in Y_1$ and $x_2 = 0 \in Y_2$. 
We identify a function $f$ on $X$ with a function $\tilde{f}$ 
on $Y_1 \sqcup Y_2$ with $\tilde{f}(z_1) = \tilde{f}(z_2)$. 
We define a bilinear form $\sE$ on 
$
\bbra{ 
  f \: ; \: 
  \tilde{f}\mbox{ is smooth on }
  Y_1 \sqcup Y_2 
} 
$ 
by
\[
\sE (f , f ) 
 = 
\frac12 
\abra{
\int_{Y_1} | \tilde{f}'(x) |^2 dx 
 + 
\int_{Y_2} | \tilde{f}'(y) |^2 dy 
}.
\]
We can easily verify that $\sE$ is closable in $L^2(X, dx)$ and 
its closure defines a Dirichlet form on $X$.
Thus we can define the corresponding diffusion process 
$( \{ Z(t) \}_{t \ge 0}, \bbra{ \P_x }_{x \in X} )$~(see \cite{FOT}).
By using identity maps 
$\iota_1 \: :\: Y_1 \to Y_2$ and $\iota_2 \: :\: Y_2 \to Y_1$, 
we define $R\: : \: X \to X$ by $R x = \iota_i(x)$ if $x \in Y_i$.
By definition of $R$ and $\P_{x_i}$($i=1,2$), (A1) and (A2) holds. 
In this case, $H = \{ z_1 \} = \{ z_2 \}$.
Let us define a map $\h \: : \: Y_2 \to Y_2$ 
by $\h(x) = 1-x$.
We define a Markovian coupling $\tilde{\P} \in \sC ( \P_{x_1} , \P_{x_2} )$ 
as the law of 
$(Z_1, Z_2)$ where $Z_1$ is a copy of $(Z , \P_{x_1})$ and 
\begin{equation}
  Z_2(t) = 
  \begin{cases}
    \h \circ R ( Z_1(t) ) & \mbox{if $t < \tau$} ,
    \\
    Z_1(t) & \mbox{if $t \ge \tau$} .
  \end{cases}
\end{equation}
Then clearly $\tilde{\P} \neq \P_M$ and 
\[
\tilde{\P} \cbra{ T (Z_1 , Z_2 ) > t }
 =
\P_{x_1} \cbra{ \tau > t }
 =
\P_M \cbra{ T (Z_1 , Z_2 ) > t }
\]
for each $t>0$. Thus $\tilde{\P}$ is also a maximal Markovian coupling.
Note that $\h \circ R$ also satisfies (A1) and (A2) 
instead of $R$ in this case.
\eexam 
\bexam{tree}{\bf (Fig.3)}
Next example is a tree. 
The space $X$, given in Fig.3, 
is a union of nine copies of the unit interval $[0,1]$ 
with some identification of these endpoints. 
$X$ is naturally regarded as a metric space. 
As in \Exam{eight}, we can construct a canonical Dirichlet form 
and the corresponding diffusion process on $X$.
Let $x_1=p_{11}$ and $x_2 = p_{22}$. There is an isometry  
$R : X \to X$ so that 
$R(p_{11}) = p_{22}$, $R(p_{12})=p_{21}$, $R(p_1)=p_2$ and 
$R$ fixes all other endpoints.
Then (A1) and (A2) holds. 
Let $\h$ be an isometry so that $\h (p_{21}) = p_{22}$ and 
$\h$ fixes all other endpoints. 
We define a Markovian coupling $\tilde{\P} \in \sC ( \P_{x_1} , \P_{x_2} )$ 
as the law of 
$(Z_1, Z_2)$ where $Z_1$ is a copy of $(Z , \P_{x_1})$ and 
\begin{equation}
  Z_2(t) = 
  \begin{cases}
    R ( Z_1 (t) ) 
     & 
    \mbox{if $t < \tau_{\{ p_1 \}}$},
    \\
    \h \circ R ( Z_1(t) ) 
     & 
    \mbox{if $\tau_{\{p_1\}} \le t < \tau_{\{ p_0\} }$},
    \\
    Z_1(t) 
     & 
    \mbox{if $t \ge \tau_{\{ p_0 \}}$} ,
  \end{cases}
\end{equation} 
where $\tau_{ \{x\} }$ is the first hitting time to $x$. 
Then clearly $\tilde{\P} \neq \P_M$ and 
\[
\tilde{\P} \cbra{ T (Z_1 , Z_2 ) > t }
=
\P_{x_1} \cbra{ \tau_{ \{ p_0 \} } >t }
=
\P_M \cbra{ T (Z_1 , Z_2 ) > t }
\]
for each $t>0$. 
Thus $\tilde{\P}$ is also a maximal Markovian coupling.
Different from \Exam{eight}, 
this example essentially has only one reflection structure.
\eexam

These examples reveal that 
maximal Markovian coupling may not be unique 
if the underlying space is more singular than 
Riemannian manifolds or Alexandrov spaces. 
One characteristic property which is common to those examples 
is the existence of branching geodesics. 
But, in general, 
non-branching property of geodesics is not necessary 
for the uniqueness of maximal Markovian coupling. 
To see this fact, we consider 
the Brownian motion on $2$-dimensional Sierpinski gasket.
\begin{figure}
\begin{center}
\begin{pspicture}(-7,-3)(7,4)
\scalebox{2}[.4]{
  \psarc[linewidth=.8pt]{<->}(0,-4.6){.2}{120}{60}
}
\rput(0.6,-2.1){$R$}
\rput(1.4,3.2){\rotatebox{90}{\scalebox{1.2}[5.9]{$\}$}}}
\rput(1.3,3.5){$X_2$}
\rput(-1.4,3.2){\rotatebox{90}{\scalebox{1.2}[5.9]{$\}$}}}
\rput(-1.3,3.5){$X_1$}
\rput(-0.3,3.8){$\hat{H}$}
\rput(-0.4,2.9){$p_3$}
\rput(-2.9,-1.5){$p_1$}
\rput(2.9,-1.5){$p_2$}
\put(1.4,0.8){$\Psi_2 (p_3)$}
\rput(3.5,-0.3){\rnode{C}{$\Psi_2(p_1)$}}
\rput(0.9,1.6){\rnode{A}{}}
\rput(5.3,2.3){
  \rnode{B}{ 
    \begin{minipage}{.28\textwidth}
      a geodesic joining $p_3$ and 
      \\ 
      $\Psi_2 \circ \Psi_1 \circ \Psi_3 \circ \Psi_2(p_3)$
    \end{minipage}
  }
}
\ncarc{<-}{A}{B}
\SpecialCoor \degrees[12]
\scalebox{1.5}{
  \pssierpinski
  \psline(2;3)(1;1)
  \psline(1;1)(.5;11)
  \rput(1;11){
    \scalebox{0.5}{ 
      \rput(1;7){
        \scalebox{0.5}{ 
          \rput(1;3){
            \scalebox{0.5}{
              \psline[linewidth=8pt](2;3)(1;1)
              \psdot[dotsize=16pt](1;1)
            }
          }
        }
      }
    }
  }
  \psset{origin={0;0}}
  \psdot[dotsize=2pt](2;3)
  \psdot[dotsize=2pt](2;7)
  \psdot[dotsize=2pt](2;11)
  \psdot[dotsize=2pt](1;9)
  \psdot[dotsize=2pt](1;1)
  \psline[linestyle=dashed,linewidth=.8pt](2.5;3)(1.5;9)
}
\rput(1.5;9){\rnode{D}{}}
\ncarc{->}{C}{D}
\end{pspicture}
\\
Fig.4
\end{center}
\end{figure}

Take three points $p_1 , p_2 , p_3 \in \R^2$ with 
$|p_i - p_j| =1$ for all $i \neq j$. 
Let us define a contraction map $\Psi_i \: : \: \R^2 \to \R^2$ 
for $i=1,2,3$ given by $\Psi_i(x) = ( x - p_i )/2 + p_i$.
Obviously, $p_i$ is the unique fixed point of $\Psi_i$.
The Sierpinski gasket is a unique compact set in $\R^2$ 
satisfying $X = \bigcup_{i=1}^3 \Psi_i(X)$~(see Fig.4). 
For detailed properties of the Sierpinski gasket, see 
\cite{Kigami} for instance. 
We set $V_0 = \{ p_1, p_2 , p_3 \}$ and 
$V_n = \bigcup_{i=1}^3 \Psi_i (V_{n-1} )$. 
The Brownian motion 
$( \{ Z(t) \}_{t \ge 0} , \{ \P_x \}_{x \in X} )$ on $X$ 
is given 
by a suitable scaling limit of 
a continuous time random walk on $V_n$ as $n\to\infty$ 
(see \cite{Barlow-Perkins,Lindstrom}). 
There is a reflection $\hat{R}$ on $\R^2$ 
so that $\hat{R}(p_1) = p_2$. 
We denote the fixed points of $\hat{R}$ by $\hat{H}$.
The map $\hat{R}$ naturally induces a reflection $R$ 
on $X$ so that its fixed points $H$ coincides 
with $X\cap \hat{H}$. 
Moreover, 
$X$ and $( \{ Z(t) \}_{t \ge 0}, \{ \P_x \}_{x\in X} )$ 
fulfills (A1) and (A2) for $x_1 = p_1$ and $x_2 = p_2$. 
As shown in \cite{Kigami}, 
there is a unique distance $d$ on $X$, 
called shortest path metric, such that it satisfies 
\begin{enumerate}
\item 
$(X,d)$ becomes a geodesic metric space,
\item 
$d( p_i , p_j ) = 1$ for each $i \neq j$,
\item 
$d( z_1 , z_2 ) = 2 d ( \Psi_i (z_1) , \Psi_i(z_2))$ 
for $z_1, z_2 \in X$ and $i=1,2,3$.
\end{enumerate}
\bthm{gasket}
Let $X$ be the Sierpinski gasket as defined above.
Then \eq{LDP_Mirror} holds. 
\ethm 
\bpf 
Let $p_t(x, y)$ be the transition density of the Brownian motion. 
Then, the main theorem of \cite{Kumag_LDP} asserts that, 
for each $u > 0$ and $x, y \in X$, 
\begin{equation}
  \label{eq:LDP_gasket}
  - \lim_{n\to\infty} \abra{ \abra{ \frac{2}{5} }^n u }^{\frac{1}{d_w-1}}
  \log p_{(2/5)^n u} (x, y) 
  = 
  d(x, y)^{d_w / (d_w - 1 )} F \abra{ \frac{u}{d(x,y)} } , 
\end{equation}
where $d_w>2$ is the walk dimension of the Sierpinski gasket and 
$F$ is an implicitly determined, non-constant, positive continuous function 
on $(0,\infty)$. 
For our aim, we need a refined observation on $F$. 
By the definition of $F$ in \cite{Kumag_LDP}, 
\begin{equation} 
F(v) = v^{1/ ( d_w - 1 )} 
\sup_{s>0} \bbra{ K(s) - vs } ,
\end{equation}
for some positive, concave and real analytic function 
$K(s)$ on $(0, \infty)$. 
Thus, 
\[
d( x, y )^{d_w / (d_w - 1 )} F \abra{ \frac{u}{d(x,y)} } 
=
u^{1 / ( d_w - 1 )} 
\sup_{s>0} \bbra{ d( x , y ) K(s) - us } . 
\]
Since $F$ is continuous on $(0,\infty)$, 
there is $s_v \in ( 0 ,\infty )$ 
for each $v \in (0,\infty)$ 
so that 
$ K(s_v) - vs_v = \sup_{s>0} \bbra{ K(s) - v s } $
holds.
Indeed, if there exists a sequence $\{ s_n \}_{n \in \N}$ 
with $\lim_{n\to\infty} s_n = \infty$ so that  
\[
\lim_{s_n\to \infty} \abra{ K(s_n) - vs_n } 
= \sup_{s >0} \bbra{ K (s) - vs },
\]
then $F(v') = \infty$ for every $v' < v$.  
These observations imply that, for $0 < a < b$,
\begin{align*}
\sup_{s>0} 
\bbra{
   a K(s) - u s 
}
& =  
a K( s_{u/a} ) - u s_{u/a}
\\
& <
b K( s_{u/a} ) - u s_{u/a}
\\
& \le 
\sup_{s>0} \bbra{ b K(s) - u s }.
\end{align*}
It means that the right hand side in \eq{LDP_gasket} 
is strictly increasing with respect to $d(x,y)$. 
Thus, 
the same argument as given in \Thm{mfd} yields that 
$x \in X_1$ and $y \in X_2$ with \eq{x>y} and \eq{x<y} satisfies 
\begin{equation}\label{eq:equidist}
d(x,z) = d(y,z) \mbox{ for all } z\in H .
\end{equation}
To complete the proof, we show that \eq{equidist} implies $x = Ry$.
It suffices to show that $w=w'$ holds 
when $w, w' \in X_2$ with $d(w ,z) = d(w', z)$ for $z=p_3$ or 
$z=\Psi_2 (p_1)$. 
In this case, we have 
\begin{equation*}
\begin{array}{rcl}
w\in \Psi_2(X)
& \Leftrightarrow 
& d(w, p_3) \ge 1/2 , d(w, \Psi_2(p_1) ) \le 1/2 ,
\\
w\in \Psi_3(X)
& \Leftrightarrow 
& d(w, p_3) \le 1/2 , d(w, \Psi_2(p_1) ) \ge 1/2 .
\end{array} 
\end{equation*} 
Thus, $w ,w'\in \Psi_2(X)$ or $w, w' \in \Psi_3(X)$.
In particular, $w=w'=\Psi_2(p_3)$
if and only if $d(w,p_3)=d(w,\Psi_2(p_1))$. 
Now we assume $w\in \Psi_2(X) \setminus \Psi_3(X)$. 
To see the argument below, we easily find that 
the same argument also works 
for the case $w, w' \in \Psi_3(X)\setminus \Psi_2(X)$. 
Since $(X,d)$ is a geodesic space, 
$d(w , p_3) = d( w, \Psi_2(p_3)) +1/2$ and therefore 
$d(w, \Psi_2(p_3) )=d ( w' , \Psi_2 (p_3) )$.
Then we have 
\begin{equation*} 
\begin{array}{rcl}
w\in \Psi_2 \circ \Psi_1(X)
& \Leftrightarrow 
& d(w, \Psi_2(p_1)) \le 1/4 , d(w, \Psi_2(p_3) ) \ge 1/4 ,
\\
w\in \Psi_2 \circ \Psi_2(X)
& \Leftrightarrow 
& d(w, \Psi_2(p_1)) \ge 1/4 , d(w, \Psi_2(p_3) ) \ge 1/4 ,
\\
w\in \Psi_2 \circ \Psi_3(X)
& \Leftrightarrow 
& d(w, \Psi_2(p_1)) \ge 1/4 , d(w, \Psi_2(p_3) ) \le 1/4 .
\end{array}
\end{equation*}
Thus $w,w' \in \Psi_2(\Psi_i(X))$ for some $i \in \{ 1,2,3\}$. 
In particular, $w=w'$ when $w\in V_2$. 
Since we have 
\begin{align*}
d(w, \Psi_2 (p_1) )
& = 
d ( w , \Psi_2 \circ \Psi_i(p_1) ) 
 + 
d( \Psi_2 \circ \Psi_i(p_1) , \Psi_2 (p_1) ) , 
\\
d(w, \Psi_2 (p_3) )
& = 
d ( w , \Psi_2 \circ \Psi_i(p_3) ) 
 + 
d( \Psi_2 \circ \Psi_i(p_3) , \Psi_2 (p_3) ) 
\end{align*}
when $w \in \Psi_2 \circ \Psi_i(X)$, 
the same argument as above works by replacing 
$\Psi_2(p_1)$, $\Psi_2(p_3)$ and $\Psi_2(X)$ 
by $\Psi_2\circ \Psi_i(p_1)$, $\Psi_2\circ \Psi_i(p_3)$ and 
$\Psi_2\circ \Psi_i(X)$ respectively.
When $w \in \bigcup_{n\in \N} V_n$, 
such a recursive argument ends in a finite step with resulting $w=w'$. 
When $w \notin \bigcup_{n \in \N} V_n$, 
we obtain a sequence $\{ i_n \}_{n \in \N}$ 
with $i_n \in \{ 1,2,3 \}$ 
so that 
$w, w' \in \Psi_{i_1}\circ \Psi_{i_2} 
\circ \cdots\circ \Psi_{i_n} (X)$ for each $n \in \N$. 
Since $\bigcap_{n \in \N} 
\Psi_{i_1}\circ \Psi_{i_2} 
\circ \cdots\circ \Psi_{i_n} (X)$ is just one point, 
$w=w'$ follows.
\epf
\brem{Sierpinski_Varadhan}
As shown in the above proof, 
\eq{LDP_gasket} means that \eq{Varadhan2} holds 
with $\ro(s) = s^{1/(d_w-1)}$, $s_n = u (2/5)^n$ and 
$\Psi(v) = v^{d_w/(d_w-1)} F(u/v)$. 
Thus the Sierpinski gasket satisfies 
all assumption in \Thm{mfd} 
except for being non-branching. 
For example, we consider 
two minimal geodesics $\gm_1$ and $\gm_2$. 
$\gm_1$ joins $p_3$ and $p_2$. 
$\gm_2$ joins $p_3$ and $\Psi_2 (p_1)$ via $\Psi_2(p_3)$. 
Then both of $\gm_1$ and $\gm_2$ contains 
the minimal geodesic 
joining $p_3$ and $\Psi_2(p_3)$. 
Thus \Thm{gasket} is not a direct consequence of \Thm{mfd}. 
\erem
\section{Kendall-Cranston coupling}
\label{sec:KC_couple}

Let $X$ be a $d$-dimensional complete Riemannian manifold 
and $( \{ Z(t) \}_{t \ge 0} , \{ \P_x \}_{x \in X} )$ 
the Brownian motion on it.
In this framework, we construct a Kendall-Cranston coupling 
following the argument due to von Renesse \cite{Renes_poly}. 
As we will see, his argument is based on approximation 
by coupled geodesic random walks. 

Let $D(X)= \bbra{ (x,x) \in X \times X \; ;\; x \in X }$. 
For each $(x,y) \in X \times X \setminus D(X)$, 
we choose a minimal geodesic $\gm_{xy}:[0,1] \to X$ of 
constant speed with $\gm_{xy}(0) =x$ and $\gm_{xy}(1)=y$. 
Let $H_{xy}$ be the hyperplane in $T_y X$ of codimension 1 
which is perpendicular to $\dot{\gm}_{xy}(1)$ and $0 \in H_{xy}$. 
For each $v \in T_x X$, take a parallel translation 
by Levi-Civita connection along $\gm_{xy}$ to $T_y X$ and 
reflect the resulting vector with respect to $H_{xy}$. 
In this way, we obtain a new vector $w \in T_y X$. 
We define a map $m_{xy} \: : \: T_x X \to T_y X$ 
by $m_{xy}v = w$. 
Clearly $m_{xy}$ is isometry. 
Take a measurable section $\ph \: : \: X \to \sO(X)$ 
to the orthonormal frame bundle $\sO(X)$. 
Let us define maps $\Ph_i \: : \: X \times X \to \sO(X)$ 
for $i=1,2$ satisfying 
\begin{align*}
\Ph_1(x,y) & \in \sO_x(X) 
& 
x,y & \in X \times X,
\\ 
\Ph_2(x,y) & \in \sO_y(X) 
& 
x,y & \in X \times X, 
\\ 
\Ph_2(x,y)u & = m_{xy} \Ph_1(x,y) u
& 
(x,y) & \in X \times X \setminus D(X), u \in \R^d,
\\ 
\Ph_1(x,x) & = \Ph_2(x,x) = \ph(x)
& 
x & \in X.
\end{align*}
We can choose $\gm_{xy}$ so that 
$(x,y) \mapsto \gm_{xy}$ is measurable
as a map from $X\times X \setminus D(X)$ to $C^1([0,1]\to X)$ 
and $\gm_{xy}$ is symmetric, i.e.~$\gm_{xy}(t) = \gm_{yx}(1-t)$.  
Also we can choose $\Ph_i$ to be measurable for $i=1,2$. 
Take a sequence of random variables $\{ \xi_n\}_{n\in \N}$ 
uniformly distributed on $d$-dimensional unit disk. 
Let us define a coupled geodesic random walk 
$Z^\ep (n) = \abra{ Z^\ep_1 (n) , Z^\ep_2 (n) }$ on $X\times X$ 
with step size $\ep>0$ and starting point $(x,y) \in X\times X$ 
inductively by 
\begin{align*}
Z^\ep (0) 
& = 
(x,y) ,
\\
Z^\ep (n+1) 
& = 
\abra{
  \exp_{Z^\ep_1 (n)} \abra{ \ep \sqrt{d+2} \Ph_1( Z^\ep (n) ) \xi_{n+1} }
  ,
  \exp_{Z^\ep_2 (n)} \abra{ \ep \sqrt{d+2} \Ph_2( Z^\ep (n) ) \xi_{n+1} }
}.
\end{align*}
Let $\tau_\lm (t)$ be the Poisson process with intensity $\lm>0$ 
independent of $\{ \xi_n \}_{n\in \N}$. 
Then the sequence of processes 
$\{ Z^{k^{-1/2}} ( \tau_k ( t ) ) \}_{ k \in \N}$ is tight 
in the Skorokhod path space $D([0,\infty) \to X \times X)$ 
and $Z^{k^{-1/2}}_i ( \tau_k ( t ) )$ weakly converges to 
the Brownian motion on $X$ as $k\to \infty$ for $i=1,2$. 
Let $\tilde{Z}(t)= ( \tilde{Z}_1 (t), \tilde{Z}_2(t) )$ be 
a (subsequential) limit of 
$\{ Z^{k^{-1/2}} ( \tau_k ( t ) ) \}_{ k \in \N}$. 
Let $\sg$ be the first hitting time of $\tilde{Z}$ to $D(X)$. 
We set $Z(t)$ by 
\begin{equation*} 
Z(t) = 
\begin{cases}
\tilde{Z}(t) 
& \mbox{if $t < \sg$,} 
\\
( \tilde{Z}_1(t), \tilde{Z}_1(t) )
& \mbox{if $t \ge \sg$.} 
\end{cases}
\end{equation*}
We call $Z(t)$ a Kendall-Cranston coupling. 
This is indeed a coupling of two Brownian motions 
starting at $x$ and $y$ respectively. 
Our choice of $\xi_n$ is a bit different from 
that in \cite{Renes_poly}, where $\xi_n$ is uniformly distributed 
on the unit sphere. 
But it does not matter since the same argument works. 

\bthm{KC}
Assume (A1) and (A2) for $x_1, x_2 \in X$. 
Then a Kendall-Cranston coupling of $\P_{x_1}$ and $\P_{x_2}$ 
is the mirror coupling defined by $R$. 
In particular, the Kendall-Cranston coupling is 
unique in the sense that it is independent of the 
choice of subsequences of approximating geodesic random walks. 
As the result, the Kendall-Cranston coupling is 
the unique maximal Markovian coupling of $\P_{x_1}$ and $\P_{x_2}$.
\ethm

\bpf
As shown in the proof of \Thm{mfd}, $R$ is an isometry on $X$. 
For $x \in X$, set $y = Rx$. 
We claim 
\begin{equation}\label{eq:ref_geod}
dR ( u )  = m_{xy} u \quad \mbox{for $u \in T_x X$}.
\end{equation}  
In order to complete the proof, 
it suffices to show \eq{ref_geod}. 
Indeed, 
since  
the equality $R ( \exp_z (w) ) = \exp_{Rz} ( dR (w) )$ holds 
for $z \in X$ and $w \in T_z X$, 
\eq{ref_geod} implies 
\begin{equation}\label{eq:ref}
Z^\ep_2 (n) = R Z^\ep_1 (n) \quad \mbox{for $n \le T(Z_1^\ep , Z_2^\ep )$}. 
\end{equation} 
Note that the coupled geodesic random walks never meet under \eq{ref}.  
That is, 
\begin{equation}\label{eq:infinite_coupling}
\P [ T (Z_1^\ep , Z_2^\ep) = \infty ] =1
\end{equation} 
for each $\ep >0$. 
This fact is shown as follows: 
by \eq{ref}, 
$Z_1^\ep ( T(Z_1^\ep , Z_2^\ep ) ) \in H$ must hold 
if $T(Z_1^\ep, Z_2^\ep) < \infty$.  
Let $\nu_{z,\ep}$ be the law of $\exp_z ( \ep \xi_1 )$. 
Then $\nu_{z,\ep}(H)=0$ for each $z \in X$ and $\ep>0$ 
since $H$ is a submanifold of codimension 1
as mentioned in \Rem{totally_geodesic} 
It easily implies \eq{infinite_coupling}. 
Once we obtain \eq{ref} and \eq{infinite_coupling}, 
the central limit theorem for geodesic random walks 
yields that the full sequence of 
$\{ Z^{k^{-1/2}} ( \tau_k ( t ) ) \}_{ k \in \N}$ 
weakly converges to the image of the Brownian motion 
by the map $z \mapsto (z, Rz)$ as $k \to \infty$. 
Thus a Kendall-Cranston coupling is unique and 
identical to the mirror coupling. 

Set $\gm = \gm_{xy}$ for simplicity. 
First we show 
\begin{equation}\label{eq:inversion}
R ( \gm (t) ) = \gm (1-t).
\end{equation} 
By the symmetric choice of $\gm_{xy}$, 
we may assume 
$l:=\inf_{z\in H} d(x,z) \le \inf_{z \in H} d(y,z)$ 
without loss of generality. 
By (A2), $l < \infty$ holds. 
Take $t_0 \in [0,1]$ so that $d(x,\gm(t_0)) = l$. 
Let $\gm_1: [0, 2 t_0] \to X$ be a curve joining $x$ and $y$ given by 
\begin{equation*}
\gm_1 (s)  = 
\begin{cases}
\gm(s)
& 
s \in [0, t_0], 
\\
R ( \gm(2t_0 - s) )
& 
s \in (t_0, 2t_0].
\end{cases}
\end{equation*}
Then the length of $\gm_1$ equals $2l$ 
and 
the minimality of $\gm$ implies 
$2l =d(x,y)$ and $t_0=1/2$. 
Moreover, $\dot{\gm}_1 (1/2) = \dot{\gm}(1/2)$ must hold 
and 
therefore $\gm_1 = \gm$. 
It proves \eq{inversion}. 
Note that the above discussion implies 
$\gm (t) \in X_1$ for $t \in [0,1/2)$ 
and $\gm(t)\in X_2$ for $t \in (1/2 , 1]$. 

Next we show \eq{ref_geod} in the case $u = \dot{\gm}(0)$. 
It easily follows from \eq{inversion}, that is, 
\begin{equation} \label{eq:paratangent}
d R ( u ) 
 = 
d R ( \dot{\gm} (0) ) 
 = 
- \dot{\gm} (1) 
 = 
m_{xy} u .
\end{equation} 
Finally we prove \eq{ref_geod} for $u \perp \dot{\gm}(0)$. 
Let $/\!\!/_{s,t} \: T_{\gm(s)} X \to T_{\gm(t)} X$ be 
the parallel translation along $\gm |_{[s,t]}$. 
It suffices to show that 
\begin{equation}\label{eq:pararef} 
/\!\!/_{1,1/2} \circ d R \circ /\!\!/_{1/2,0}(h) = h 
\end{equation}
for each $h \in T_{\gm (1/2)} X$ 
with $ h \perp \dot{\gm} (1/2)$. 
Indeed, once we prove it, 
\begin{align*}
dR ( u )
 = 
d R ( /\!\!/_{1/2,0} \circ /\!\!/_{0,1/2} ( u ) )
 = 
/\!\!/_{1/2,1} \circ /\!\!/_{0,1/2} ( u ) 
 = 
/\!\!/_{0,1} ( u )  
 =
m_{xy} u .
\end{align*}
Now we show \eq{pararef}. Take $\ep >0$ so that 
the exponential map $\exp_{\gm(1/2)} \: : \: T_{\gm (1/2)} X \to X$ 
is diffeomorphic on $2\ep$-ball centered at $0 \in T_{\gm(1/2)}X$.  
We may assume that $| h | = \ep$. 
Let $ h' = \exp_{\gm(1/2)}^{-1} ( \gm (1/2-\ep / d ( \gm(0) , \gm(1) ) ) )$. 
Note that $/\!\!/_{1,1/2} \circ d R \circ /\!\!/_{1/2,0} (h') = - h'$.
We consider a curve $c : [0,1] \to X$ given by 
 $c(t) = \exp ( \cos \pi t \, h' + \sin \pi t \, h )$. 
Since $c(0) \in X_1$ and $c(1) \in X_2$, (A2) yields that 
$c$ intersects $H$. 
By the choice of $h'$, $R ( c(t) ) \neq c(t)$ if $t \neq 1/2$. 
Hence $c(1/2) \in H$. 
It implies \eq{pararef} and therefore completes the proof.
\epf 

\noindent 
\textit{Acknowledgment.}
The problem treated in section \ref{sec:KC_couple} was 
suggested by K.-Th.~Sturm. 
The author would like to thank him for his advice.

\bibliographystyle{amsplain}
\bibliography{refs}

\end{document}